\documentclass{article}
\usepackage{xy}
\usepackage{epsfig}
\usepackage[T2A]{fontenc} 
\usepackage[cp1251]{inputenc} 
\usepackage[russian]{babel}
\usepackage{amsfonts,amssymb,amsmath}
\usepackage{hyperref} 
\usepackage{graphics} 
\usepackage{multido, xcolor} 

\newtheorem{theorem}{\bf\large Òåîðåìà}[section]

\newtheorem{proposition}[theorem]{\large\bf Ïðåäëîæåíèå}

\newtheorem{conjecture}[theorem]{\large\bf Ãèïîòåçà}

\makeindex

\def \F {{\mathbb F}}

\def\ph{\varphi}

\def\is{\buildrel \sim   \over \rightarrow}

\newcommand{\M}{{{\rm M}}}

\newcommand{\ord}{{\rm ord}}
\newcommand{\SL}{{\rm SL}}

\newcommand{\N}{{\mathbb N}}
\newcommand{\Q}{{\mathbb Q}}
\newcommand{\R}{{\mathbb R}}
\newcommand{\PP}{{\mathbb P}}
\newcommand{\Z}{{\mathbb Z}}

\newcommand{\al}{\alpha}

\newcommand{\e}{\varepsilon}
\newcommand{\ep}{\varepsilon}

\font\teneusm=eusm10 \font\seveneusm=eusm7 
\font\fiveeusm=eusm5 
\newfam\eusmfam 
\textfont\eusmfam=\teneusm 
\scriptfont\eusmfam=\seveneusm 
\scriptscriptfont\eusmfam=\fiveeusm

\def\mat #1,#2,#3,#4,{\left({#1\atop #3}{#2\atop #4}\right)}
\def\bra#1,{{\left\lbrace {#1}\right\rbrace}}
\def\ph{\varphi}

\def \M{{\rm M}}

\def\is{\buildrel \sim   \over \rightarrow}

\def\l1{\langle}

\newcommand{\B}{\left(\begin{array}{cc}}
\newcommand{\E}{\end{array}\right)}

\def \fns{{${}^{*)}$}}

\usepackage{euscript}
\let\scr=\EuScript

\let\mathcal=\scr           
\def\ang#1,{{\left\langle {#1}\right\rangle}} 
\newcommand{\ds}{\displaystyle}

\def \F {{\mathbb F}}

\def\ph{\varphi}

\def\is{\buildrel \sim   \over \rightarrow}

\font\teneusm=eusm10 \font\seveneusm=eusm7 
\font\fiveeusm=eusm5 
\newfam\eusmfam 
\textfont\eusmfam=\teneusm 
\scriptfont\eusmfam=\seveneusm 
\scriptscriptfont\eusmfam=\fiveeusm

\font\tengothic=eufm10
\font\sevengothic=eufm7
\font\fivegothic=eufm5
\newfam\Gothic
\textfont\Gothic=\tengothic
\scriptfont\Gothic=\sevengothic
\scriptscriptfont\Gothic=\fivegothic

\def\mat #1,#2,#3,#4,{\left({#1\atop #3}{#2\atop #4}\right)}
\def\bra#1,{{\left\lbrace {#1}\right\rbrace}}
\def\ph{\varphi}

\def\w{\omega}

\def \Re {{\rm Re}}

\def \Card {{\rm Card}\,}

\def \M{{\rm M}}

\def\is{\buildrel \sim   \over \rightarrow}

\def\l1{\langle}

\let\scr=\EuScript

\let\mathcal=\scr           

\def\vin{{ {\tiny \mid }  
\kern-7.29pt 
\bigcup }}

\def\ang#1,{{\left\langle {#1}\right\rangle}}

\newcommand{\cross}{\times}
\normalsize

\newcommand{\CC}{\mathbb C}

\newcounter{ncours}{\setcounter{ncours} {1}}

\def\SL {\mathop{\rm SL}\nolimits}

\font\tenrus=wncyr10
\font\sevenrus=wncyr7
\newfam\Rus
\textfont\Rus=\tenrus
\scriptfont\Rus=\sevenrus

\def\scr{\fam\Rus\sevenrus} 

\font\tenbrus=wncyb10
\font\eightbrus=wncyb8
\newfam\Brus
\textfont\Brus=\tenbrus
\scriptfont\Brus=\eightbrus

\font\tenirus=wncyi10
\font\eightirus=wncyi8
\newfam\Irus
\textfont\Irus=\tenirus
\scriptfont\Irus=\eightirus

\font\tenrus=wncyr10
\font\sevenrus=wncyr7
\newfam\Rus
\textfont\Rus=\tenrus
\scriptfont\Rus=\sevenrus

\def\scr{\fam\Rus\sevenrus} 
\input cyracc.def
\def\i1{\accent'044i}
\def\I1{\accent'044I}
\def\e1{\accent'040e}
\def\l1{l{}p1}

\title{Ëîêàëüíûå è ãëîáàëüíûå ìåòîäû â àðèôìåòèêå
 }

\author 
{À.À.Ïàí÷èøêèí}
\date{}

\begin{document}
\maketitle


%
%

\begin{abstract}
Ïóñòü $p$ -- ïðîñòîå ÷èñëî. 
Îáñóæäàþòñÿ ìåòîäû ðåøåíèÿ ñðàâíåíèé ïî ìîäóëþ $p^n$ ñ ïîìîùüþ  
$p$-àäè÷åñêèõ ÷èñåë, àíàëîãè÷íûå ìåòîäàì ðåøåíèÿ óðàâíåíèé â äåéñòâèòåëüíûõ  
÷èñëàõ.
Ïðèâåäåíû ïðèìåðû ñâÿçè ëîêàëüíûõ è ãëîáàëüíûõ ìåòîäîâ 
ðåøåíèÿ â àðèôìåòèêå,
à òàêæå ïðèìåðû êîìïüþòåðíûõ âû÷èñëåíèé ñ $p$-àäè÷åñêèìè ÷èñëàìè è àëãåáðàè÷åñêèìè êðèâûìè.

\end{abstract}
  \tableofcontents
\section{Ââåäåíèå}

Ñòàòüÿ îñíîâàíà íà ìàòåðèàëàõ ëåêöèé àâòîðà â Èíñòèòóòå Ôóðüå (Ãðåíîáëü, Ôðàíöèÿ), â Ýêîëü Íîðìàëü (Ëèîí, Ôðàíöèÿ),
à òàêæå íà ìàòåðèàëàõ ñïåöêóðñîâ íà ìåõ-ìàòå ÌÃÓ â 1979-1991 è â 2001.

Â ñòàòüå îáñóæäàþòñÿ ñëåäóþùèå òåìû:

\begin{enumerate}
\item[1)] 
{$p$-àäè÷åñêèå ÷èñëà è ñðàâíåíèÿ.} 
Äèîôàíòîâû ñèñòåìû ëèíåéíûõ óðàâíåíèé è ñðàâíåíèé
\item[2)] 
{Ïðèíöèï ÌèíêîâñêîãîÕàññå äëÿ êâàäðàòè÷íûõ ôîðì}
\item[3)] 
{Ñèìâîë Ãèëüáåðòà è åãî âû÷èñëåíèå} 
\item[4)]
{Êóáè÷åñêèå óðàâíåíèÿ è ýëëèïòè÷åñêèå êðèâûå} 
\item[5)] 
{Îò ñðàâíåíèé ê ðàöèîíàëüíûì òî÷êàì: ãèïîòåçà Á¸ð÷à è Ñóèííåð\-òîíà--Äàéåðà}
\end{enumerate}


\section{$p$-àäè÷åñêèå ÷èñëà è ñðàâíåíèÿ} 
Èäåÿ ðàñøèðåíèÿ ïîëÿ $\Q$
 â òåîðèè ÷èñåë
âñòðå÷àåòñÿ â ðàçëè÷íûõ âàðèàíòàõ. 
Íàïðèìåð, âëîæåíèå, $\Q \subset \R$
 ÷àñòî äàåò ïîëåçíûå íåîáõîäèìûå óñëîâèÿ ñóùåñòâîâàíèÿ ðåøåíèé äèîôàíòî\-âûõ
óðàâíåíèé íàä $\Q$ è íàä $\Z$.
Âàæíîå ñâîéñòâî ïîëÿ $\R$ 
iåãî ïîëíîòà: ëþáàÿ
ôóíäàìåíòàëüíàÿ ïîñëåäîâàòåëüíîñòü (ïîñëåäîâàòåëüíîñòü Êîøè):
 $\bra {\al_n},_{n=1}^\infty$ â $\R$
 èìååò ïðåäåë $\al$.
 Ôóíäàìåíòàëü\-íîñòü îçíà÷àåò, ÷òî äëÿ ïðîèçâîëüíîãî $\ep > 0$ ìàëû àáñîëþòíûå
  âåëè÷èíû $|\al_n - \al_m|<\ep$ äëÿ âñåõ $n$ è $m$
  áîëüøèõ íåêîòîðîãî íàòóðàëüíîãî ÷èñëà $N = N(\ep)$.
Êðîìå òîãî, âñå ýëåìåíòû $\R$ ÿâëÿþòñÿ ïðåäåëàìè ôóíäàìåíòàëüíûõ ïîñëåäîâàòåëüíîñòåé
 $\bra {\al_n},_{n=1}^\infty$ ñ $\al_n \in \Q$.

Àíàëîãè÷íàÿ êîíñòðóêöèÿ ñóùåñòâóåò è äëÿ âñåõ
 $p$-àäè÷åñêèõ íîðìèðîâàíèé
 $|\cdot|_p$ ïîëÿ $\Q$: 
\begin{align*}&
|\cdot|_p:\Q \to \R_{\ge  0}=\bra {x\in \R\ |x\ \ge 0}, \cr
&|a/b|_p=p^{\ord_pb-\ord_pa},\ \ |0|_p=0,
\end{align*}
 ãäå $\ord_pa$ íàèâûñøàÿ ñòåïåíü ÷èñëà $p$ äåëÿùàÿ öåëîå ÷èñëî $a$.
Ýòà îáùàÿ êîíñòðóêöèÿ
«ïðèñîåäèíåíèÿ ïðåäåëîâ ôóíäàìåíòàëüíûõ ïîñëåäîâàòåëüíîñòåé» îòíîñèòåëüíî
 íåêîòîðîãî íîðìèðîâàíèÿ  $|\cdot|$ ïîëÿ $k$  
 íàçûâàåòñÿ ïîïîëíåíèåì.
Â ðåçóëüòàòå ïîëó÷àåòñÿ ïîëå  $\hat k$,
ñ íîðìèðî\-âà\-íèåì, òàêæå îáîçíà÷àåìîì $|\cdot|$
 ïðè÷åì ïîëå $\hat k$ -- ïîëíîå, à $k$ îäíîçíà÷íî âêëàäûâàåòñÿ â
 $\hat k$ â êà÷åñòâå âñþäó
ïëîòíîãî ïîäïîëÿ ñ ñîõðàíåíèåì íîðìèðîâàíèÿ, ñì. 
\cite{Borevich Z.I. Shafarevich I.R. (1985)}, 
\cite{Koblitz N. (1980)}.

Ñîãëàñíî òåîðåìå Îñòðîâñêîãî, âñå íîðìèðîâàíèÿ ïîëÿ $\Q$
 ñâîäÿòñÿ ëèáî ê àáñîëþòíîé âåëè÷èíå, ëèáî ê $p$-àäè÷åñêîìó íîðìèðîâàíèþ
 $|\cdot|_p$ (ñ òî÷íîñòüþ äî ýêûèûàëåíòíîñòè).
Ïîýòîìó âñå ïîïîëíåíèÿ ïîëÿ $\Q$
 ýòî ëèáî ïîëå äåéñòâèòåëüíûõ ÷èñåë, $\R$,  ëèáî ïîëÿ $\Q_p$ 
 $p$-àäè÷åñêèõ ÷èñåë.
Èñïîëüçîâàíèå âñåâîçìîæíûõ âëîæåíèé
 $\Q \hookrightarrow \R$ è $\Q \hookrightarrow \Q_p$
 ($p$--ïðîñòîå ÷èñëî)
 ÷àñòî çíà÷èòåëüíî óïðîùàåò ñèòóàöèþ â àðèôìåòè÷åñêèõ çàäà÷àõ.
Çàìå÷àòåëüíûé ïðèìåð äàåò {\it òåîðåìà  ÌèíêîâñêîãîÕàññå} 
 (ñì.\cite{Borevich Z.I. Shafarevich I.R. (1985)}, ãëàâà 1).
Óðàâíåíèå
\begin{align}\label{ii23.1}
 Q(x_1, x_2, \dots, x_n)=0, 
\end{align}
 çàäàííîå êâàäðàòè÷íîé ôîðìîé $Q(x_1, x_2, \dots, 
x_n)=\sum_{i,j}a_{ij}x_ix_j$,
 $a_{ij}\in \Q$, èìååò
íåòðèâèàëüíîå ðåøåíèå â ðàöèîíàëüíûõ ÷èñëàõ â òîì è òîëüêî â òîì ñëó÷àå, êîãäà îíî íåòðèâèàëüíî ðàçðåøèìî
 íàä $\R$ è íàä $\Q_p$ äëÿ âñåõ ïðîñòûõ ÷èñåë $p$.
Äëÿ íàõîæäåíèÿ ðåøåíèé óðàâíåíèé íàä $\Q_p$ ìîæíî ýôôåêòèâíî
ïðèìåíÿòü ïðèåìû, âçÿòûå ïî àíàëîãèè èç àíàëèçà íàä $\R$,
 òàêèå, êàê «ìåòîä êàñàòåëüíûõ Íüþòîíà»
 ({\it ``Newton - Raphson algorithm''}), 
 êîòîðûé â $p$--àäè÷åñêîì ñëó÷àå èçâåñòåí êàê  {\it ëåììà Ãåíçåëÿ,
 (Hensel's lemma)}.\index{Hensel's lemma}

Íàèáîëåå ïðîñòûì ñïîñîáîì ìîæíî ââåñòè $p$-àäè÷åñêèå ÷èñëà êàê âûðàæåíèÿ âèäà
\begin{align}\label{ii23.2}
 \al = a_m p^m + a_{m+1} p^{m+1} + \dots, 
\end{align}
ãäå $a_i\in \{0, 1, \dots. p-1\}$--öèôðû (ïî îñíîâàíèþ $p$, à 
$m\in \Z$. Óäîáíî çàïèñûâàòü $\al$ â âèäå ïîñëåäîâàòåëüíîñòè öèôð, áåñêîíå÷íîé âëåâî:
$$
\al = \begin{cases}\cdots a_{m+1}a_m{\buildrel m-1 \ {\rm zeros } \over 
{\overbrace {000\dots 0}}_{(p)}},&\mbox{if }m\ge 0, \cr
\cdots a_1a_0.a_{-1} \cdots a_m{}_{(p)}, &\mbox{if }m<0.
\end{cases}
$$
Ýòè âûðàæåíèÿ îáðàçóþò ïîëå, â êîòîðîì îïåðàöèè âûïîëíÿþòñÿ òàê æå,
êàê äëÿ íàòóðàëüíûõ ÷èñåë $n=a_0+a_1p+\dots a_rp^r$, çàïèñàííûõ ïî îñíîâàíèþ $p$.
Ñëåäîâàòåëüíî, â ýòîì ïîëå ëåæàò íàòóðàëüíûå, à ïîòîìó è âñå
ðàöèîíàëüíûå ÷èñëà. Íàïðèìåð,
$$
 -1
 =
 \frac {p-1}{ 1-p}
 =
 (p-1) + (p-1)p+(p-1)p^2+\cdots
 =
 \cdots (p-1) (p-1)_{(p)};
$$
$$
 \frac {-a_0}{ p-1}
 =
 a_0+a_0p+a_0p^2+\cdots
 =
 \cdots a_0a_0a_0{}_{(p)}.
$$
Åñëè $n\in \N$, òî âûðàæåíèå äëÿ $-n=n\cdot(-1)$ âèäà 
(\ref{ii23.2}) ïîëó÷àåòñÿ, åñëè ïåðåìíîæèòü óêàçàííûå âûðàæåíèÿ äëÿ $n$ è äëÿ $-1$.
Âîîáùå, åñëè $\al \in \Q$ òî çàïèøåì $\al = c-\frac a b $, ãäå $a,c\in 
\Z$, $b\in N$, $0\le a<b$, ò.å. $a/b$ ïðàâèëüíàÿ äðîáü. Òîãäà ïî ýëåìåíòàðíîé òåîðåìå Ýéëåðà, $p^{\ph(b)}-1=bu$, $u\in \N$.
Ïîýòîìó
$$
-\frac a b  = \frac {au}{1 - p^{\ph(b)}}  ,
$$
è $au<bu=p^r-1$, $r=\ph (b)$.   Òåïåðü ìû âèäèì, ÷òî çàïèñü ïî îñíîâàíèþ
$p$ ÷èñëà $au$ èìååò âèä $a_{r-1}\cdots a_0{}_{(p)}$,  
ñëåäîâàòåëüíî, âûðàæåíèå (\ref{ii23.2}) 
äëÿ ÷èñëà $\al$ ïîëó÷àåòñÿ êàê ñóììà âûðàæåíèÿ äëÿ $c\in \N$ è   
$$
-\frac a b  = \cdots a_0{\buildrel r \ {\rm digits }\over {\overbrace 
{a_{r-1}\cdots a_0}}}{\buildrel r\ {\rm digits}\over {\overbrace 
{a_{r-1}\cdots a_0}}}{}_{(p)}.
$$
Íàïðèìåð, äëÿ $p=5$ èìååì  
$$
\frac 9 7  = 2-\frac 5 7  = 2 + \frac {5\cdot {2232}}  {1-5^6} \ \ c=2\ 
a=5,\ b=7,
$$
ïðè÷åì
$$
2232=32412_{(5)} = 3\cdot 5^4 + 2\cdot 5^3 + 4\cdot 5^2 + 1\cdot 5 + 
2,
$$
ïîýòîìó
$$
\frac 9 7  = \cdots  {\overbrace {324120}}{\overbrace 
{324120}}324122_{(5)}.
$$
Íåòðóäíî ïðîâåðèòü, ÷òî ïîïîëíåíèå ïîëÿ $\Q$ îòíîñèòåëüíî $p$--àäè÷åñêîé
ìåòðèêè $|\cdot|_p$ îòîæäåñòâëÿåòñÿ ñ ïîëåì «$p$-àäè÷åñêèõ ðàçëîæåíèé» âèäà
  (\ref{ii23.2}). 
Ïðè ýòîì $|\al|_p=p^m$   ãäå â âûðàæåíèè
 (\ref{ii23.2}) äëÿ $\al$ èìååì $a_m\not=0$ (ñì. \cite{Koblitz N. (1980)}).
 
 $p$-àäè÷åñêèå ðàçëîæåíèÿ ìîæíî ðàññìàòðèâàòü êàê àíàëîãè ðàçëîæåíèÿ ôóíêöèè $f$ ïåðåìåííîé $x$ 
 â îêðåñòíîñòè òî÷êè  $a$ ïî ñòåïåíÿì  $(x-a)$  , ïðè÷¸ì   $p$  ÿâëÿåòñÿ àíàëîãîì  $(x-a)$:
 \subsubsection*{Âû÷èñëåíèå ñ PARI-GP}
\begin{verbatim}
gp > forprime(p=2,163,print("p="p,",""9/7="9/7+O(p^6)))
p=2,9/7=1 + 2 + 2^2 + 2^3 + 2^5 + O(2^6)
p=3,9/7=3^2 + 3^3 + 2*3^5 + O(3^6)
p=5,9/7=2 + 2*5 + 5^2 + 4*5^3 + 2*5^4 + 3*5^5 + O(5^6)
p=7,9/7=2*7^-1 + 1 + O(7^6)
p=11,9/7=6 + 9*11 + 7*11^2 + 4*11^3 + 9*11^4 + 7*11^5 + O(11^6)
p=13,9/7=5 + 9*13 + 3*13^2 + 9*13^3 + 3*13^4 + 9*13^5 + O(13^6)
p=17,9/7=11 + 14*17 + 4*17^2 + 7*17^3 + 2*17^4 + 12*17^5 + O(17^6)
p=19,9/7=4 + 8*19 + 5*19^2 + 16*19^3 + 10*19^4 + 13*19^5 + O(19^6)
p=23,9/7=21 + 9*23 + 16*23^2 + 19*23^3 + 9*23^4 + 16*23^5 + O(23^6)
p=29,9/7=22 + 20*29 + 20*29^2 + 20*29^3 + 20*29^4 + 20*29^5 + O(29^6)
p=31,9/7=19 + 26*31 + 8*31^2 + 13*31^3 + 4*31^4 + 22*31^5 + O(31^6)
p=37,9/7=33 + 15*37 + 26*37^2 + 31*37^3 + 15*37^4 + 26*37^5 + O(37^6)
p=41,9/7=13 + 29*41 + 11*41^2 + 29*41^3 + 11*41^4 + 29*41^5 + O(41^6)
p=43,9/7=32 + 30*43 + 30*43^2 + 30*43^3 + 30*43^4 + 30*43^5 + O(43^6)
p=47,9/7=8 + 20*47 + 13*47^2 + 40*47^3 + 26*47^4 + 33*47^5 + O(47^6)
p=53,9/7=24 + 45*53 + 37*53^2 + 22*53^3 + 45*53^4 + 37*53^5 + O(53^6)
p=59,9/7=35 + 50*59 + 16*59^2 + 25*59^3 + 8*59^4 + 42*59^5 + O(59^6)
p=61,9/7=10 + 26*61 + 17*61^2 + 52*61^3 + 34*61^4 + 43*61^5 + O(61^6)
p=67,9/7=30 + 57*67 + 47*67^2 + 28*67^3 + 57*67^4 + 47*67^5 + O(67^6)
p=71,9/7=52 + 50*71 + 50*71^2 + 50*71^3 + 50*71^4 + 50*71^5 + O(71^6)
p=73,9/7=43 + 62*73 + 20*73^2 + 31*73^3 + 10*73^4 + 52*73^5 + O(73^6)
p=79,9/7=69 + 33*79 + 56*79^2 + 67*79^3 + 33*79^4 + 56*79^5 + O(79^6)
p=83,9/7=25 + 59*83 + 23*83^2 + 59*83^3 + 23*83^4 + 59*83^5 + O(83^6)
p=89,9/7=14 + 38*89 + 25*89^2 + 76*89^3 + 50*89^4 + 63*89^5 + O(89^6)
p=97,9/7=29 + 69*97 + 27*97^2 + 69*97^3 + 27*97^4 + 69*97^5 + O(97^6)
p=101,9/7=59 + 86*101 + 28*101^2 + 43*101^3 + 14*101^4 + 72*101^5 + O(101^6)
p=103,9/7=16 + 44*103 + 29*103^2 + 88*103^3 + 58*103^4 + 73*103^5 + O(103^6)
p=107,9/7=93 + 45*107 + 76*107^2 + 91*107^3 + 45*107^4 + 76*107^5 + O(107^6)
p=109,9/7=48 + 93*109 + 77*109^2 + 46*109^3 + 93*109^4 + 77*109^5 + O(109^6)
p=113,9/7=82 + 80*113 + 80*113^2 + 80*113^3 + 80*113^4 + 80*113^5 + O(113^6)
p=127,9/7=92 + 90*127 + 90*127^2 + 90*127^3 + 90*127^4 + 90*127^5 + O(127^6)
p=131,9/7=20 + 56*131 + 37*131^2 + 112*131^3 + 74*131^4 + 93*131^5 + O(131^6)
p=137,9/7=60 + 117*137 + 97*137^2 + 58*137^3 + 117*137^4 + 97*137^5 + O(137^6)
p=139,9/7=41 + 99*139 + 39*139^2 + 99*139^3 + 39*139^4 + 99*139^5 + O(139^6)
p=149,9/7=129 + 63*149 + 106*149^2 + 127*149^3 + 63*149^4 + 106*149^5 + O(149^6)

p=151,9/7=66 + 129*151 + 107*151^2 + 64*151^3 + 129*151^4 + 107*151^5 + O(151^6)

p=157,9/7=91 + 134*157 + 44*157^2 + 67*157^3 + 22*157^4 + 112*157^5 + O(157^6)
p=163,9/7=141 + 69*163 + 116*163^2 + 139*163^3 + 69*163^4 + 116*163^5 + O(163^6)
\end{verbatim}
Ëþáîïûòíî ñðàâíèòü ðàçëîæåíèÿ (\ref{ii23.2}) «áåñêîíå÷íûå âëåâî», ñ
 ðàçëîæåíèÿìè äåéñòâèòåëüíûõ ÷èñåë $\al \in \R$, «áåñêîíå÷íûìè âïðàâî»: 
\begin{align}\label{pinfty}
\al = a_ma_{m-1} \cdots a_0.a_{-1}\cdots = a_m10^m + a_{m-1}10^{m-1} 
+ \cdots a_0 + a_{-1}10^{-1} + \cdots,
\end{align}
ãäå $a_i \in \bra {0,1,\cdots, 9},$ öèôðû, à $a_m\not=0$. 
Ðàçëîæåíèÿ òàêîãî òèïà ïî ëþáîìó
íàòóðàëüíîìó îñíîâàíèþ ïðèâîäÿò ê îäíîìó è òîìó æå ïîëþ $\R$, ïðè ýòîì
îíè íåîäíîçíà÷íû, ê ïðèìåðó, $2.000\cdots = 
1.999\cdots$. 
Ðàçëîæåíèÿ  (\ref{pinfty}) ìîæíî ðàññìàòðèâàòü êàê àíàëîãè ðàçëîæåíèÿ 
 â îêðåñòíîñòè òî÷êè  $p=\infty$, ïðè÷¸ì   $p=\infty$ ÿâëÿåòñÿ àíàëîãîì  $x^{-1}$.

Ðàçëîæåíèÿ (\ref{ii23.2}) â $p$--àäè÷åñêîì ñëó÷àå âñåãäà îäíîçíà÷íî îïðåäåëåíû, ÷òî ñîçäàåò äîïîëíèòåëüíûå âû÷èñëèòåëüíûå
óäîáñòâà.

Ïîëå  $\Q_p$ ÿâëÿåòñÿ {\it ïîëíûì ìåòðè÷åñêèì ïðîñòðàíñòâîì} 
 ñ òîïîëîãèåé, îïðåäåëåííîé ñèñòåìîé «îòêðûòûõ äèñêîâ» âèäà:
$$
U_a(r) = \bra {x\ | \ |x-a|<r}, \ \ \ (x,\ a\in \Q_p,\ r>0)
$$
(èëè «çàìêíóòûõ äèñêîâ» $D_a(r) = \bra {x\ | \ |x-a|\le r},$).
Ïðè ýòîì è $U_a(r)$ è $D_a(r)$ ÿâëÿþòñÿ îòêðûòî-çàìêíóòûìè ìíîæåñòâàìè ñ òîïîëîãè÷åñêîé òî÷êè çðåíèÿ $\Q_p$.

Âàæíîå òîïîëîãè÷åñêîå ñâîéñòâî ïîëÿ 
 $\Q_p$ 
 åãî {\it  ëîêàëüíàÿ êîìïàêòíîñòü}: âñå äèñêè êîíå÷íîãî ðàäèóñà êîìïàêòíû. Â ýòîì ïðîùå âñåãî
 óáåäèòüñÿ íà ÿçûêå ïîñëåäîâàòåëüíîñòåé, ïîêàçàâ, ÷òî êàæäàÿ ïîñëåäîâàòåëüíîñòü $\bra {\al_n},_{n=1}^\infty$ ýëåìåíòîâ äèñêà $\al_n \in D_a(r)$
 èìååò â ýòîì æå äèñêå ïðåäåëüíóþ òî÷êó. Ýòà ïðåäåëüíàÿ òî÷êà ëåãêî èùåòñÿ ñ ïîìîùüþ
 {\it $p$---àäè÷åñêèõ
öèôð}
 (\ref{ii23.2}) ïîñëåäîâàòåëüíî, ñïðàâà íàëåâî, è èñïîëüçóåòñÿ òîò ôàêò, ÷òî
ó âñåõ ýëåìåíòîâ $\al_n \in D_a(r)$ ÷èñëî çíàêîâ «ïîñëå çàïÿòîé» îãðàíè÷åíî ôèêñèðîâàííûì ÷èñëîì. 
Â ÷àñòíîñòè, äèñê
$$
 \Z_p 
 =
 D_0(1)
 =
 \bra {x\ |\ |x|_p\le 1},
 =
 \bra { x=a_0 + a_1 p + a_2 p^2 + \cdots },
$$
 -- ýòî êîìïàêòíîå òîïîëîãè÷åñêîå êîëüöî, 
 ýëåìåíòû êîòîðîãî íàçûâàþòñÿ öåëûìè
 $p$--àäè÷åñêèìè ÷èñëàìè, ïðè ýòîì
 $\Z_p$ ñîâïàäàåò ñ çàìûêàíèåì ìíîæåñòâà îáû÷íûõ öåëûõ ÷èñåë $\Z$ â $\Q_p$.
Êîëüöî $\Z_p$ ÿâëÿåòñÿ ëîêàëüíûì,
ò. å. èìååò åäèíñòâåííûé ìàêñèìàëüíûé èäåàë $p\Z_p = U_0(1)$
 ñ ïîëåì âû÷åòîâ $\Z_p/p\Z_p = \F_p$.
Ìíîæåñòâî îáðàòèìûõ ýëåìåíòîâ åäèíèö êîëüöà $\Z_p$ -- ýòî
$$
 \Z_p^\times
 =
 \Z_p \backslash p\Z_p
 =
 \bra {x\ |\ |x|_p = 1},
 =
 \bra {x = a_0 + a_1 p + a_2 p^2+ \cdots\ |\ a_0\not=0},.
$$
Äëÿ êàæäîãî ýëåìåíòà $x\in \Z_p$ îïðåäåëåí åãî ïðåäñòàâèòåëü Òåéõìþëëåðà
 $$
 \displaystyle \w (x)  = \lim_{n\to \infty}x^{p^n}.
 $$
(ïðåäåë âñåãäà ñóùåñòâóåò è óäîâëåòâîðÿåò óðàâíåíèþ: $\w(x)^p = \w 
(x)$, è ñïðàâåäëèâî ñðàâíåíèå $\w (x) \equiv x \bmod\ p$. Íàïðèìåð, äëÿ $p=5$  èìååì
\begin{align*}&\w(1) = 1;\cr &
\w (2) = 2 + 1\cdot 5 + 2\cdot 5^2 + 1\cdot 5^3 + 3\cdot 5^4 \cdots; 
\cr &
\w (3) = 3 + 3\cdot 5 + 2\cdot 5^2 + 3\cdot 5^3 + 1\cdot 5^4 + 
\cdots; \cr &
\w(4) = 4 + 4\cdot5 + 4\cdot 5^2 + 4\cdot 5^3 + 4\cdot 5^4 + \cdots = 
-1; \cr &
\w (5) = 0.
\end{align*}

Êîëüöî $\Z_p$ ìîæíî îïèñàòü òàêæå êàê ïðîåêòèâíûé ïðåäåë 
 $$
 \lim_{\longleftarrow\atop n} \Z / p^n\Z
 $$
êîëåö  $A_n = \Z/p^n\Z$
îòíîñèòåëüíî ãîìîìîðôèçìîâ ðåäóêöèè  ïî ìîäóëþ $p^{n-1}$
 $\ph_n: A_n \to A_{n-1}$.
Ïîñëåäîâàòåëüíîñòü 
\begin{align}\label{EQii23.3}
 \cdots
 {\buildrel \ph_{n+1} \over \longrightarrow }
 A_n
 {\buildrel \ph_{n} \over \longrightarrow }
 A_{n-1}
 {\buildrel \ph_{n-1} \over \longrightarrow }
 \cdots
 {\buildrel \ph_3 \over \longrightarrow }
 A_2
 {\buildrel \ph_2 \over \longrightarrow }
 A_1
\end{align}
 îáðàçóåò ïðîåêòèâíóþ ñèñòåìó, çàíóìåðîâàííóþ öåëûìè ÷èñëàìè
  $n\ge 1$.
Ïðîåêòèâíûé ïðåäåë ñèñòåìû (\ref{EQii23.3})ýòî êîëüöî
 $$
 \lim_{\longleftarrow\atop n}A_n
 $$
 ñî ñëåäóþùèì
óíèâåðñàëüíûì ñâîéñòâîì: îäíîçíà÷íî îïðåäåëåíû òàêèå ãîìîìîðôèçìû
ïðîåêöèè
 $$
 \pi_n : \lim_{\longleftarrow\atop n}A_n \to A_n,
 $$
÷òî äëÿ ïðîèçâîëüíîãî êîëüöà $B$ è ñèñòåìû ãîìîìîðôèçìîâ
 $\psi _n: B \to A_n$ ñîãëàñîâàííûõ äðóã ñ äðóãîì óñëîâèåì:
 $\psi_{n-1} = \ph_n \circ \psi_n$ for $n\ge 2$,
 ñóùåñòâóåò åäèíñòâåííûé ãîìîìîðôèçì  $\psi: B \to A$
 äëÿ êîòîðîãî $\psi_n = \pi_n \psi$ 
(ñì. \cite{Koblitz N. (1980)}, 
\cite{Serre J.--P. (1970)}).
Äëÿ êîëüöà $A$ ïîñòðîåíèå ãîìîìîðôèçìîâ ïðîåêöèè
è ïðîâåðêà óíèâåðñàëüíîãî ñâîéñòâà íåïîñðåäñòâåííî âûòåêàþò
èç çàïèñè åãî ýëåìåíòîâ ñ ïîìîùüþ «öèôð» 
(\ref{ii23.2}). 

Àíàëîãè÷íî,
 $$
 \Z_p^\times
 =
 \lim_{\longleftarrow\atop n}(\Z/p^n\Z)^\times,
 $$
 (ïðîåêòèâíûé ïðåäåë ãðóïï). Äëÿ îïèñàíèÿ ãðóïïû
 $\Q_p^{\times}$ ïîëîæèì $\nu =1$ åñëè $p>2$ è $\nu =2$ åñëè $p=2$,
è îïðåäåëèì
 $$
 U = U_p
 =
 \bra {x\in \Z_p |x\equiv 1\,\bmod\,p^{\nu}},.
 $$
Òîãäà $U\is \Z_p$
 (èçîìîðôèçì ìóëüòèïëèêàòèâíîé è àääèòèâíîé ãðóïï)$U_p$
 è $\Z_p$. Äëÿ
ïîñòðîåíèÿ ýòîãî èçîìîðôèçìà çàìåòèì, ÷òî 
$$
 U \is \lim_{\longleftarrow\atop n}U/U^{p^n}
$$ 
è îïðåäåëèì èçîìîðôèçìû êîíå÷íûõ ãðóïï
$$
 \al_{p^n}:U/U^{p^n}\is \Z/p^n\Z,
$$
ïîëîæèâ
\begin{align}\label{ii23.4}
 \al_{p^n}((1+p^{\nu})^a)
 =
 a\,\bmod\,p^n\ \ (a\in \Z).
\end{align}
Ïðîñòàÿ ïðîâåðêà ïîêàçûâàåò, ÷òî îòîáðàæåíèÿ (\ref{ii23.4}) êîððåêòíî îïðåäåëåíû è ÿâëÿþòñÿ èçîìîð\-ôèçìàìè. 
Òàêèì îáðàçîì, ãðóïïà $U$ýòî òîïîëîãè÷åñêàÿ öèêëè÷åñêàÿ ãðóïïà, â êà÷åñòâå îáðàçóþùåé êîòîðîé ìîæíî âçÿòü 
$1+p^{\nu}$. Äðóãîå äîêàçàòåëüñòâî ñëåäóåò èç ñâîéñòâ ôóíêöèè, îïðåäåëåííîé
ñòåïåííûì ðÿäîì
$$
\log (1+x) = \sum_{n=1}^\infty (-1)^{n+1}\frac {x^n}  n ,
$$
\index{Logarithm!p-adic --@Logarithm!$p$-adic --}
\index{p-adic!-- logarithm@$p$-adic!-- Logarithm}
êîòîðàÿ çàäà¸ò èçîìîðôèçì èç $U$ íà $p\Z_p$

Ïðîâåðÿåòñÿ, ÷òî ñïðàâåäëèâû ðàçëîæåíèÿ
\begin{align}\label{ii23.5}
 \Q_p^\times
 =
 p^{\Z}\times \Z_p^{\times},\ \
 \Z_p^{\times}
 \cong
 (\Z/p^{\nu}\Z)^{\times} \times U.
\end{align}

\subsection{Ïðèëîæåíèÿ $p$--àäè÷åñêèõ ÷èñåë ê ðåøåíèþ ñðàâíåíèé.}\label{ii23.2.}
Âîçíèêíîâåíèå $p$--àäè÷åñêèõ ÷èñåë â ðàáîòàõ Ãåíçåëÿ áûëî òåñíî ñâÿçàíî
ñ ïðîáëåìîé ðåøåíèÿ ñðàâíåíèé ïî ìîäóëþ $p^n$, 
à ïðèìåíåíèå èõ ê òåîðèè êâàäðàòè÷íûõ ôîðì åãî ó÷åíèêîì Õàññå ïðèâåëî ê ýëåãàíòíîé ôîðìóëèðîâêå òåîðèè êâàäðàòè÷íûõ ôîðì íàä ðàöèîíàëüíûìè ÷èñëàìè, íå
èñïîëüçóþùåé ðàññìîòðåíèé â êîëüöàõ âû÷åòîâ âèäà  $\Z/p^n\Z$,
ðàáîòàòü ñ êîòîðûìè çàòðóä\-íè\-òåëü\-íî èç-çà íàëè÷èÿ äåëèòåëåé íóëÿ â $\Z/p^n\Z$. 
Èç ïðåäñòàâëåíèÿ êîëüöà $\Z_p$ â âèäå ïðîåêòèâíîãî ïðåäåëà $\Z/p^n\Z$, 
$$
\lim_{\longleftarrow\atop n}\Z/p^n\Z
$$ 
 âûòåêàåò, ÷òî
åñëè $f(x_1, \dots, x_n) \in \Z_p [x_1, \dots, x_n]$, òî ñðàâíåíèÿ
$$
f(x_1, \dots, x_n) \equiv 0 (\bmod\ p^n)
$$
 ðàçðåøèìû ïðè ëþáîì $n\ge 1$ òîãäà è òîëüêî òîãäà, êîãäà óðàâíåíèå
$$
f(x_1, \dots, x_n) = 0
$$
 ðàçðåøèìî â öåëûõ $p$--àäè÷åñêèõ ÷èñëàõ.
Ýòè ðåøåíèÿ â $\Z_p$ ìîæíî íàõîäèòü
ñ ïîìîùüþ $p$-àäè÷åñêîãî âàðèàíòà ìåòîäà êàñàòåëüíûõ Íüþòîíà 
 ({\it ``Newton - Raphson algorithm''}).

\begin{theorem}[ëåììà Ãåíçåëÿ]\index{Hensel's Lemma}
Ïóñòü $f(x)\in \Z_p[x]$ -- ìíîãî÷ëåí îäíîé ïåðåìåííîé $x$,  $f^\prime (x) \in \Z_p [x]$ 
 -- åãî ôîðìàëüíàÿ ïðîèçâîäíàÿ, è äëÿ
íåêîòîðîãî
$\al_0 \in \Z_p$ âûïîëíåíî íà÷àëüíîå óñëîâèå
\begin{align}\label{ii23.6}
 |f(\al_0)/ f^\prime (\al_0)^2|_p < 1 
\end{align}

Òîãäà ñóùåñòâóåò åäèíñòâåííîå òàêîå $\al\in \Z_p$, ÷òî
$$
 f(\al) = 0,\ \ |\al - \al_0|<1.
$$
\end{theorem}

\

Äîêàçàòåëüñòâî ïðîâîäèòñÿ ñ ïîìîùüþ ðàññìîòðåíèÿ ïîñëåäîâàòåëüíîñòè:
$$
\al_n = \al_{n-1} - \frac {f(\al_{n-1})}{ f^\prime (\al_{n-1})} .
$$
Ñ ó÷åòîì ôîðìàëüíîãî ðàçëîæåíèÿ Òåéëîðà ìíîãî÷ëåíà $f(x)$
â òî÷êå  
$x=\al_{n-1}$
 ïðîâåðÿåòñÿ, ÷òî ïîñëåäîâàòåëüíîñòü ôóíäàìåíòàëüíà,
à å¸ ïðåäåë $\al$ îáëàäàåò âñåìè íåîáõîäèìûìè ñâîéñò\-âà\-ìè (ñì. 
\cite{Borevich Z.I. Shafarevich I.R. (1985)}, 
\cite{Serre J.--P. (1970)}).

Íàïðèìåð, åñëè $f(x) = x^{p-1}-1$,
 òî ëþáîå $\al_0 \in \bra {1, 2, \dots, p-1}, $
 óäîâëåòâîðÿåò óñëîâèþ $|f(\al_0)|_p<1$,
â òî âðåìÿ êàê
 $f^\prime(\al_0) = (p-1)\al_0^{p-2} \not \equiv 0\bmod\ p$,
 ïîýòîìó íà÷àëüíîå óñëîâèå (\ref{ii23.6}) âûïîëíåíî.
Êîðåíü $\al$ ñîâïàäàåò ïðè
ýòîì ñ åäèíñòâåííûì ïðåäñòàâèòåëåì Òåéõìþëëåðà ÷èñëà $\al_0$: $\al = \w (\al_0)$.

Îïèñàííûé ìåòîä ïðèìåíèì è ê ìíîãî÷ëåíàì ìíîãèõ ïåðåìåííûõ, íî óæå
áåç ñîõðàíåíèÿ åäèíñòâåííîñòè íàõîäèìîãî ðåøåíèÿ, (ñì. 
 \cite{Borevich Z.I. Shafarevich I.R. (1985)}, 
 \cite{Koblitz N. (1980)}, 
\cite{Serre J.--P. (1970)}). 

Åùå îäíî ïðèëîæåíèå ëåììû Ãåíçåëÿ ñâÿçàíî ñ îïèñàíèåì êâàäðàòîâ
ïîëÿ $\Q_p$: äëÿ ïðîèçâîëü\-íîãî ýëåìåíòà
$$
\al = p^m\cdot v\ \in\ \Q_p\ \ (m\in \Z,\ v\in\Z_p^\times),
$$
 ñâîéñòâî $\al$ áûòü êâàäðàòîì â $\Q_p$ ðàâíîñèëüíî òîìó, ÷òî
\medskip
\begin{description}
\item{a)}
åñëè $p>2$, òî $m\in 2\Z$, à $\overline v = v \bmod\ p \in 
(\Z/p\Z)^\times{}^2$ (ò.å. $\left( \frac {\overline v}  p  \right) = 1$, 
ãäå $\left( \frac {\overline v} p \right)$ ñèìâîë Ëåæàíäðà
\item{b)}
 åñëè $p=2$, òî $m\in 2\Z$, à $v\equiv 1 \bmod\ 8$.
\end{description}
\medskip
Ðàçðåøèìîñòü óðàâíåíèÿ $x^2 = \al$ â $\Q_p$ â óñëîâèÿõ à) è á) âûâîäèòñÿ
èç ëåììû Ãåíçåëÿ, à íåîáõîäèìîñòü èõ âûòåêàåò èç áîëåå òðèâèàëüíûõ
ðàññìîòðåíèé ïî ìîäóëþ $p$ è ïî ìîäóëþ 8.

Êàê ñëåäñòâèå ìû ïîëó÷àåì, ÷òî ôàêòîðãðóïïà
 $\Q_p^\times / \Q_p^\times {}^2$ 
\medskip
\begin{description}
\item{a)}
 ïðè $p>2$ èçîìîðôíà  $\Z/2\Z \times \Z/2\Z$ ñ ñèñòåìîé ïðåäñòàâèòåëåé $\bra {1, p, v,  pv},$, $\left( \frac{ \overline v} 
p \right) = -1$;
\item{b)}
 ïðè $p=2$ èçîìîðôíà $\Z/2\Z \times \Z/2\Z \times \Z/2\Z $  
ñ ñèñòåìîé ïðåäñòàâèòåëåé $\bra {\pm1, \pm5, \pm 2, 
\pm 10},$.
\end{description}

\section{Äèîôàíòîâû ñèñòåìû ëèíåéíûõ óðàâíåíèé è ñðàâíåíèé}\label{i12.2.}
\index{Diophantine equation!Linear --}
\subsection{Âû÷èñëåíèÿ ñ êëàññàìè âû÷åòîâ.}\label{i11.4.}
Ñ òî÷êè çðåíèÿ àëãåáðû,
ìíîæåñòâî öåëûõ ÷èñåë ${\Bbb Z}$ ÿâëÿåòñÿ êîììóòàòèâíûì àññîöèàòèâíûì êîëüöîì ñ åäèíèöåé, ò. å. ìíîæåñòâîì ñ äâóìÿ êîììóòàòèâíûìè è àññîöèàòèâíûìè îïåðàöèÿìè (ñëîæåíèå è óìíîæåíèå), ñâÿçàííûìè äðóã ñ äðóãîì
çàêîíîì äèñòðèáóòèâíîñòè.
Ïîíÿòèå äåëèìîñòè â êîëüöàõ ñâÿçàíî ñ ïîíÿòèåì èäåàëà. 
Èäåàëîì $I$ â êîììóòàòèâíîì àññîöèàòèâíîì êîëüöå $R$ íàçûâàåòñÿ
ïîäìíîæåñòâî ñ
$RIR\subset I$.

Èäåàë âèäà $I=aR,\ a\in A$ íàçûâàåòñÿ ãëàâíûì èäåàëîì, ïîðîæäåííûì ýëåìåíòîì $(a)$. 
Òîãäà îòíîøåíèå äåëèìîñòè $a|b$ â êîëüöå $R$ ðàâíîñèëüíî
âêëþ÷åíèþ ñîîòâåòñòâóþùèõ ãëàâíûõ èäåàëîâ:
$$
 (b)\subset(a) \quad \ {\mbox{ èëè  } }\ \quad b\in(a).
$$
Â êîëüöå ${\Bbb Z}$ äåëåíèå ñ îñòàòêîì íà íàèìåíüøèé ïîëîæèòåëüíûé ýëåìåíò
â èäåàëå $I\ne 0$ ïîêàçûâàåò, ÷òî âñå èäåàëû ãëàâíûå, ò. å. âñÿêèé íåíóëåâîé èäåàë $I$ èìååò âèä $(N) = N{\Bbb Z}$ äëÿ íàòóðàëüíûõ ÷èñåë $N > 1$. 
Ïðè ýòîì
èäåàëû, ìàêñèìàëüíûå ïî âêëþ÷åíèþ, â òî÷íîñòè ñîîòâåòñòâóþò ïðîñòûì
÷èñëàì. Îñòàòêè îò äåëåíèÿ íà $N$ ïîäðàçäåëÿþò âñå öåëûå ÷èñëà íà íåïåðåñåêàþùèåñÿ êëàññû
$$
 \bar a=a+N\Z,\quad 0\le a \le N-1,
$$
ìíîæåñòâî êîòîðûõ òàêæå îáðàçóåò êîëüöî, îáîçíà÷àåìîå
$$
 \Z/N\Z=\Z/(N)=\{\bar 0,\bar 1,\dots,\overline{N-1}\},
$$
è ïèøåòñÿ $a\equiv b \ (\bmod \    N)$ âìåñòî $\bar 
a = \bar b$.
×àñòî â çàäà÷àõ òåîðèè ÷èñåë âû÷èñëåíèÿ â êîëüöå ${\Bbb Z}$
ìîæíî ñâîäèòü ê âû÷èñëåíèÿì â êîëüöå âû÷åòîâ
$\Z/N\Z$. 
Ýòî äîñòàâëÿåò ðÿä óäîáñòâ, íàïðèìåð, íà ìíîãèå ýëåìåíòû èç $\Z/N\Z$ ìîæíî äåëèòü, îñòàâàÿñü â ïðåäåëàõ
ýòîãî êîëüöà (â îòëè÷èå îò öåëûõ ÷èñåë, ãäå âñåãäà îïðåäåëåíî òîëüêî
äåëåíèå íà $\pm 1$).
Äåéñòâèòåëüíî, åñëè ÷èñëî $a$ âçàèìíî ïðîñòî ñ $N$, ò. å.
$\gcd (a,N)=1$,  êëàññ $\bar a$ îáðàòèì, òàê êàê â ýòîì ñëó÷àå 
ñóùåñòâóþò òàêèå öåëûå ÷èñëà $x,\  y$, ÷òî $ax+Ny=1$, ïîýòîìó $\bar a\cdot \bar x = \bar 1$.
Òàê ïîëó÷àþòñÿ âñå îáðàòèìûå ýëåìåíòû  êîëüöà âû÷åòîâ $\Z/ N\Z$
 êîòîðûå îáðàçóþò ãðóïïó ïî óìíîæåíèþ, îáîçíà÷àåìóþ $(\Z/ N\Z)^\times$. 
Ïîðÿäîê ýòîé ãðóïïû, îáîçíà÷àåòñÿ $\varphi (N)$
(ôóíêöèÿ Ýéëåðà). Íàçâàíèå ïðîèñõîäèò èç îáîáùåíèÿ ìàëîé òåîðåìû Ôåðìà, ïðèíàäëåæàùåãî Ýéëåðó:
\begin{align}
 {a}^{\varphi (N)}\equiv 1 (\bmod \    N)  \label{i11.4}
\end{align}
äëÿ âñåõ òàêèõ ýëåìåíòîâ $a$, ÷òî $\gcd (a,N)=1$,
ò.å. 
 $\bar {a}^{\varphi (N)}=\bar 1$ äëÿ òàêèõ ýëåìåíòîâ
  $\bar a$
â êîëüöå ${\Bbb Z}/N{\Bbb Z}$.

Äîêàçàòåëüñòâî Ýéëåðà, ïðèìåíèìîå ê ëþáîé êîíå÷íîé àáåëåâîé ãðóïïå 
ïîðÿäêà $f$, ïîêàçû\-âà\-åò, ÷òî ïîðÿäîê ëþáîãî ýëåìåíòà $a$ äåëèò $f$.
Óìíîæåíèå íà $a$ ÿâëÿåòñÿ ïåðåñòàíîâêîé ìíîæåñòâà ýëåìåíòîâ êîíå÷íîé àáåëåâîé ãðóïïû
(à äàííîì ñëó÷àå  ãðóïïû $(\Z/ N\Z)^\times$ ïîðÿäêà $f=\ph (N)$). 
Ïðîèçâåäåíèå âñåõ ýëåìåíòîâ ãðóïïû
óìíîæàåòñÿ íà $a^f$  ïðè ýòîé ïåðåñòàíîâêå. Ïîýòîìó $a^f=1$.

Åñëè ÷èñëî $N$ ðàçëîæåíî â ïðîèçâåäåíèå $N=N_1N_2\dots N_k$  ïîïàðíî âçàèìíî
 ïðîñòûõ ÷èñåë $N_i$, òî èìååòñÿ ðàçëîæåíèå
\begin{align}
 {\Z}/N{\Z} \cong {\Z}/N_1{\Z} \oplus
\dots \oplus {\Z}/N_k{\Z}. \label{i11.5}
\end{align}
â ïðÿìîå ïðîèçâåäåíèå êîëåö, ÷òî ýêâèâàëåíòíî êèòàéñêîé òåîðåìå îá
îñòàòêàõ: äëÿ ëþáûõ âû÷åòîâ $a_i \bmod \    N_i,\quad i=1,\dots,k$ 
íàéäåòñÿ òàêîå öåëîå ÷èñëî $a$, ÷òî  
$a\equiv a_i \bmod \    N_i$ äëÿ âñåõ  $i$.
Ïðàêòè÷åñêèé ïîèñê ÷èñëà $a$
ìîæíî áûñòðî îñóùåñòâèòü, ïðèìåíÿÿ ïîâòîðíî àëãîðèòì Åâêëèäà. Ïîëîæèì  $M_i=N/N_i$, 
òîãäà ÷èñëà $M_i$ and $N_i$
ïî óñëîâèþ âçàèìíî ïðîñòû, è 
ñóùåñòâóþò òàêèå
öåëûå ÷èñëà $X_i$ ÷òî $X_iM_i\equiv 1 \bmod \    N_i$. 
Ïîëîæèì òåïåðü
\begin{align}
 a=\sum_{i=1}^k a_iX_iM_i. \label{i11.6}
\end{align}
Ñëåäîâàòåëüíî, ÷èñëî $a$ èñêîìîå. 
Êðîìå òîãî, èç ðàçëîæåíèÿ (\ref{i11.5}) âûòåêàåò è ðàçëîæåíèå ìóëüòèïëèêàòèâíîé
ãðóïïû:
\begin{align}
 (\Z/N\Z)^\times
 \cong
 (\Z/N_1\Z)^\times \times \dots \times (\Z/N_k\Z)^\times,
\label{i11.7}
\end{align}
èç êîòîðîãî, â ÷àñòíîñòè, ñëåäóåò, ÷òî 
$\varphi (N)=\varphi (N_1) \dots \varphi (N_k)$. 
Ïîñêîëüêó äëÿ ïðîñòîãî ÷èñëà 
$p$ èìååì $\varphi (p^a)=p^{a-1}(p-1)$, ìû íàõîäèì $\varphi (N)$ 
èñõîäÿ èç ðàçëîæåíèÿ ÷èñëà  $N$.

Â ñïåöèàëüíîì ñëó÷àå, êîãäà $N=q$ ïðîñòîå ÷èñëî, êîëüöî âû÷åòîâ $\Z/N\Z$
ÿâëÿåòñÿ ïîëåì: â íåì îáðàòèì ëþáîé ýëåìåíò, îòëè÷íûé îò íóëÿ.

\subsection{Óðàâíåíèå $ax + by = c$}
Â ýòîì ïàðàãðàôå âñå áóêâû (êîýôôèöèåíòû è íåèçâåñòíûå â óðàâíåíèÿõ) îçíà÷àþò öåëûå ÷èñëà. 
Ìíîæåñòâî
$$
I(a,b)=\{c \quad |\quad \mbox{ óðàâíåíèå }  ax+by=c  \mbox{ ðàçðåøèìî (â öåëûõ ÷èñëàõ)}\}
$$
ÿâëÿåòñÿ èäåàëîì êîëüöà ${\Bbb Z}$ è ïîýòîìó $I(a,b)$
èìååò âèä
$d {\Bbb Z}$, ãäå $d={ÍÎÄ}(a,b)$ -- íàèáîëüøèé îáùèé äåëèòåëü. Òàêèì îáðàçîì,
óðàâíåíèå
\begin{align}
ax+by=c  \label{i12.1}
\end{align}
ðàçðåøèìî, òîëüêî åñëè $d$ äåëèò $c$. Êîíêðåòíîå ðåøåíèå íàõîäèòñÿ
 ñ ïîìîùüþ àëãîðèòìà Åâêëèäà: åñëè $X$, $Y$ ñ
$aX+bY=d$ òî ÷èñëà $x_0=eX, y_0=eY$ óäîâëåòâîðÿþò óðàâíåíèþ, ãäå $e=c/d$. 
Òåïåðü ìû ïîëó÷èëè âñå öåëî÷èñëåííûå ðåøåíèÿ:
$$
 x=x_0+(b/d)t,\quad y=y_0-(a/d)t,
$$
ãäå $t$ ïðîèçâîëüíîå öåëîå ÷èñëî.

Óðàâíåíèå (\ref{i12.1}) äàåò ïåðâûé ïðèìåð îáùåé ïðîáëåìû:
äëÿ ñèñòåìû óðàâíåíèé, çàäàííîé öåëî÷èñëåííûìè ìíîãî÷ëåíàìè
\begin{align}
 F_1(x_1,\dots ,x_n)=0,\quad\cdots ,\quad F_m(x_1,\dots ,x_n)=0
  \label{i12.2}
\end{align}
íàéòè âñå öåëî÷èñëåííûå (èëè âñå ðàöèîíàëüíûå) ðåøåíèÿ. Äëÿ óðàâíåíèÿ (\ref{i12.1})
 çàäà÷à íàõîæäåíèÿ ðàöèîíàëüíûõ ðåøåíèé òðèâèàëüíà. Åñëè
â ñèñòåìå (\ref{i12.2})
âñå óðàâíåíèÿ $F_i = 0$ ëèíåéíûå, òî è äëÿ íåå âñå ðàöèîíàëüíûå ðåøåíèÿ
 ëåãêî íàõîäÿòñÿ ïîñëåäîâàòåëüíûì èñêëþ÷åíèåì íåèçâåñòíûõ (íàïðèìåð, ïî ìåòîäó Ãàóññà).


\subsection{Ñèñòåìû ëèíåéíûõ óðàâíåíèé ñ öåëûìè êîýôôèöèåíòàìè}\label{i11.4.}
Îïèøåì îáùèé
ïðèåì íàõîæäåíèÿ âñåõ öåëî÷èñëåííûõ ðåøåíèé ñèñòåìû (öåëî÷èñëåííûõ)
ëèíåéíûõ óðàâíåíèé, çàïèñàííîé â ìàòðè÷íîé ôîðìå
\begin{align}
Ax=b, \label{i12.3}
\end{align}
ãäå
$$
A=\left( \begin{matrix}a_{11} & a_{12} &\cdots & a_{1n}\cr
\cdots &\cdots &\ddots &\cdots\cr 
a_{m1}&a_{m2}&\cdots&a_{mn}\end{matrix} \right) \in M_{m,n}(\Bbb Z),\;\;
x= \left( \begin{matrix}x_1\cr \cdots \cr x_n\cr 
\end{matrix} \right),\;\;
b= \left( \begin{matrix} b_1\cr \cdots \cr b_m \cr \end{matrix} \right).
$$
Ñ ïîìîùüþ òåîðèè ýëåìåíòàðíûõ äåëèòåëåé ìàòðèöû ýòà çàäà÷à òàêæå
ñâîäèòñÿ ê ïðèìåíåíèþ àëãîðèòìà Åâêëèäà. Ýëåìåíòàðíûì ïðåîáðàçîâàíèåì íàä ${\Bbb Z}$
 ñòðîê ìàòðèöû íàçîâåì ïðåîáðàçî\-âà\-íèå, ïðè êîòîðîì ê íåêîòîðîé ñòðîêå ïðèáàâëÿþò äðóãóþ, óìíîæåííóþ íà öåëîå ÷èñëî, à îñòàëüíûå
ñòðîêè íå ìåíÿþò. Ïðîâåðÿåòñÿ, ÷òî ïðèìåíåíèå òàêîãî ïðåîáðàçîâàíèÿ
ýêâèâàëåíòíî óìíîæåíèþ èñõîäíîé ìàòðèöû ñëåâà íà íåêîòîðóþ ìàòðèöó
 $U=E_{ij}=E+\lambda e_{ij}$
èç $\SL_m({\Bbb Z})$ (ñîîòâ. $\SL_n({\Bbb Z}))$ 
(öåëî÷èñëåííóþ ìàòðèöó ñ îïðåäåëèòåëåì, ðàâíûì 1). Àíàëîãè÷íîå ïðåîáðàçîâàíèå ñòîëáöîâ ðàâíîñèëüíî óìíîæåíèþ ìàòðèöû ñïðàâà íà $V\in \SL_m({\Bbb Z})$

Ïðèìåíåíèå íåñêîëüêèõ òàêèõ ïðåîáðàçîâàíèé ñ ýëåìåíòàðíûìè
ìàòðèöàìè ïðèâîäèò ìàòðèöó $A$ ê âèäó $UAV$, à öåëî÷èñëåííûå ðåøåíèÿ ñîîòâåòñòâóþùåé ñèñòåìû óðàâíåíèé
\begin{align} 
UAVy=Ub \label{i12.4} 
\end{align}
è èñõîäíîé ñèñòåìû (\ref{i12.3}) âçàèìíî îäíîçíà÷íî ñîîòâåòñòâóþò äðóã äðóãó
ïî ôîðìóëå $x=Vy$. Òåïåðü íàèáîëüøèé îáùèé äåëèòåëü $d_1$ ýëåìåíòîâ ìàòðèöû $A$ ìîæíî íàéòè ïîâòîðíûì ïðèìåíåíèåì àëãîðèòìà Åâêëèäà ê åå
ýëåìåíòàì $a_{i,j}$ , èñïîëüçóÿ ýëåìåíòàðíûå ïðåîáðàçîâàíèÿ ñòðîê è ñòîëáöîâ
è ïðè íåîáõîäèìîñòè ìåíÿÿ çíàê ñòðîêè òàê, ÷òî ïðåîáðàçîâàííàÿ ìàòðèöà $A^\prime$ ïðèìåò âèä
\begin{align}
D=\left( \begin{matrix} d_1&0&0&\dots &0\cr
                 0&d_2&0&\dots &0\cr
                 \cdots&\cdots &\ddots &\dots &0\cr
                 0&0&\dots &d_r &\dots \cr
                   \cdots &\cdots &\cdots &\cdots &\dots\cr
\end{matrix}
 \right) = UAV.  \label{i12.5} 
\end{align}
Òåïåðü ìû ïîëó÷àåì ðåøåíèå ïðåîáðàçîâàííîé (à ïîýòîìó è èñõîäíîé)
öåëî÷èñëåííîé ñèñòåìû ëèíåéíûõ óðàâíåíèé:
$d_iy_i=c_i,\quad c=Ub$ äëÿ $i \le r,$ $c_i = 0$ äëÿ îñòàëüíûõ $i$, ïðè ýòîì
$y_i$ ïðèíèìàþò ïðîèçâîëüíûå öåëûå
çíà÷åíèÿ. 
Êðèòåðèé ñîâìåñòíîñòè íàä  ${\Bbb Z}$ ñîñòîèò â òîì, ÷òî $d_i|c_i$ äëÿ âñåõ $i\le r$, 
 è $c_i = 0$ äëÿ îñòàëüíûõ $i$.

×èñëà $d_i$ íàçûâàþòñÿ ýëåìåíòàðíûìè äåëèòåëÿìè ìàòðèöû $A$.
Ïðîèçâåäåíèÿ $d_1\cdots d_i$ 
 ñîâïàäàþò ñ íàèáîëüøèìè îáùèìè äåëèòåëÿìè âñåõ ìèíîðîâ ïîðÿäêà $i$
ìàòðèöû $A$ è $d_i|d_{i+1}.$

Îòñþäà ñëåäóåò è òàêàÿ ôîðìóëèðîâêà êðèòåðèÿ ñîâìåñòíîñòè íàä $\Z$ ñèñòåìû (\ref{i12.3}): 
äëÿ ýòîãî íåîáõîäèìî è äîñòàòî÷íî, ÷òîáû
áûëà ðàçðåøèìà ñîîòâåòñòâóþùàÿ ñèñòåìà ñðàâíåíèé
$$
 Ax \equiv b (\bmod \ N)
$$
ïî ëþáîìó íàòóðàëüíîìó ìîäóëþ $N > 2$. 
Êðèòåðèé òàêîãî ðîäà íàçûâàåòñÿ
ïðèíöèïîì Ìèí\-êîâñ\-êîãî  Õàññå è îí ÷àñòî âñòðå÷àåòñÿ â çàäà÷àõ äèîôàíòîâîé ãåîìåòðèè.

\section{Óðàâíåíèÿ âòîðîé ñòåïåíè. }\label{i12.3.}
\subsection{Êâàäðàòè÷íûå ôîðìû è êâàäðèêè}
\index{Diophantine equation!Quadratic --}
Äëÿ äèîôàíòîâà óðàâíåíèÿ
\begin{align}
 f(x_1,x_2,\dots 
,x_n)=\sum_{i,j}^na_{ij}x_ix_j+\sum_{i=1}^nb_ix_i+c=0.
 \label{i12.6}
\end{align}
íàõîäèòü öåëî÷èñëåííûå ðåøåíèÿ çíà÷èòåëüíî òðóäíåå, ÷åì ðàöèîíàëüíûå,
õîòÿ è ýòà çàäà÷à óæå íåòðèâèàëüíà. 
Èçâåñòíûé ïðèìåð ñâÿçàí ñ ðàöèîíàëüíîé
 ïàðàìåòðèçàöèåé îêðóæíîñòè $x^2+y^2=1:$
 ïî ôîðìóëàì óíèâåðñàëüíîé ïîäñòàíîâêè
\begin{align}
 x=\frac {2t} {1+t^2} \ ,
 y=\frac {1-t^2} {1+t^2} 
 \quad (x=\cos \varphi , \
 y=\sin \varphi ,\ 
 t=\tan \left(\frac {\varphi} {2} \right)).\label{i12.7}
\end{align}
èç êîòîðîé ñëåäóåò îïèñàíèå âñåõ ïðèìèòèâíûõ ïèôàãîðåéñêèõ òðîåê
 $(X,Y,Z)$,
ò.å. íàòóðàëüíûõ ðåøåíèé $X^2+Y^2=Z^2$
ñ  ${ ÍÎÄ} (\ X ,\ Y ,\ Z\ )$ $ = 1$
ïî ôîðìóëàì: $X=2uv,\  Y=u^2-v^2,\  Z=u^2+v^2,$
ãäå $u>v>0$ âçàèìíî ïðîñòûå ÷èñëà ïðîòèâîïîëîæíîé ÷¸òíîñòè. 
Äëÿ ýòîãî íàäî â ôîðìóëàõ (\ref{i12.7}) ïîëîæèòü $t=u/v$ .  

Âîîáùå, ïðè îòûñêàíèè ðàöèîíàëüíûõ ðåøåíèé óðàâíåíèÿ (\ref{i12.6}) óäîáíî
ïåðåéòè ê êâàäðà\-òè÷\-íîé ôîðìå
\begin{align}
 &F(X_0, X_1,\cdots , X_n)
  =
 \sum_{i,j=0}^nf_{ij}X_iX_j \cr
 &  \hskip1cm =
 \sum_{i,j=1}^nf_{ij}X_iX_j+2\sum_{i=1}^nf_{i0}X_iX_0+f_{00}X_0^2,\label{i12.8}
\end{align}
ãäå $f_{ij}=f_{ji}=a_{ij}/2$ äëÿ $1\le i<j\le n$ è 
$f_{0i}=f_{i0}=b_i/2$
ñ $i=1,2,\dots ,n,\; f_{00}=c$.
Äëÿ ýòîãî íàäî çàìåíèòü íåîäíîðîäíûå êîîðäèíàòû  $x_1,\dots ,x_n$ íà îäíîðîäíûå $X_0,\dots ,X_n$ 
ïî ôîðìóëàì $X_i=x_iX_0\quad
(i=1,2,\dots ,n).$
Êâàäðàòè÷íàÿ ôîðìà $F(X)$ ÿâëÿåòñÿ îäíîðîäíûì ìíîãî÷ëåíîì âòîðîé ñòåïåíè,
 êîòîðûé óäîáíî çàïèñûâàòü â ìàòðè÷íîé ôîðìå
$$
 F(X)=X^tA_FX, \quad X^t=(X_0,X_1,\dots ,X_n),
$$
ãäå $A_F=(f_{ij})$ ìàòðèöà êîýôôèöèåíòîâ.
Åñëè ñóùåñòâóåò íåíóëåâîå ðàöèîíàëüíîå ðåøåíèå $F(X)=0$, òî ãîâîðÿò, ÷òî
 $F$ ïðåäñòàâëÿåò íóëü
íàä ïîëåì ðàöèîíàëüíûõ ÷èñåë. 
Ýòî óðàâíåíèå îïðåäåëÿåò êâàäðèêó $Q_F$,
êîòîðóþ ìû áóäåì ðàññìàòðèâàòü êàê ãèïåðïîâåðõíîñòü â êîìïëåêñ\-íîì
ïðîåêòèâíîì ïðîñòðàíñòâå ${\CC}{\PP}^n$: 
$$
 Q_F
 =
 \{(z_0:z_1:\dots :z_n)\in {\CC}{\PP}^n \quad |\quad F(z_0, z_1,\dots , 
z_n)=0 \}.
$$
Íåíóëåâîå ðàöèîíàëüíîå ðåøåíèå $F(X)=0$ îïðåäåëÿåò òî÷êó $X_0$
íà êâàäðèêå $Q_F$. 
Îñòàëüíûå ðàöèîíàëüíûå òî÷êè (ðàöèîíàëüíûå ðåøåíèÿ) ëåãêî íàéòè:
 îíè ñîâïàäàþò ñ òî÷êàìè ïåðåñå÷åíèÿ êâàäðèêè $Q_F$ ñî
âñåâîçìîæíûìè ïðÿìûìè, âûõîäÿùèìè èç $X_0$
è îïðåäåëåííûìè íàä ${\Bbb Q}$ (ò.å. â íàïðàâëåíèè âåêòîðà ñ
 ðàöèîíàëüíûìè êîîðäèíàòàìè). 
Ïðÿìàÿ ïðîõîäÿùàÿ ÷åðåç $X^0$ è $Y^0$ ñîñòîèò èç òî÷åê $uX^0+vY^0$.
Óðàâíåíèå $F(uX^0+vY^0)=0$ ñâîäèòñÿ ê
$$
uv\sum_{i=1}^n\frac {\partial F} {\partial 
X_i} (X^0)Y_i^0+v^2F(Y^0)=0.
$$
Íàäî òîëüêî, ÷òîáû òî÷êà $X^0$ íå áûëà «âåðøèíîé» íà $Q_F$, ò.å.  $\frac {\partial 
F} {\partial X_i} (X^0)\ne 0$
 õîòÿ áû äëÿ îäíîãî $i$.
Â ýòîì ñëó÷àå, äëÿ ëþáîãî  $Y^0$ íàõîäèòñÿ òî÷êà ïåðåñå÷åíèÿ  $Q_F$ ñ ýòîé ïðÿìîé:
\begin{align}
v=-u\sum_{i=1}^n\frac {\partial F} {\partial 
X_i} (X^0)Y_i^0/F(Y^0). \label{i12.9}
\end{align}
(Åñëè  $F(Y^0)=0$ òî $Y^0$ óæå  íà $Q_F$). 
Ïðèìåðîì ðàññìîòðåííîé êîíñòðóêöèè, çàïèñàííûì
â íåîäíîðîäíûõ êîîðäèíàòàõ, ÿâëÿþòñÿ
ôîðìóëû  (\ref{i12.7}): ÷òîáû íàéòè âñå ïàðû $(x, y)$
ðàöèîíàëüíûõ ÷èñåë, äëÿ êîòîðûõ $x^2+y^2=1$, ðàññìîòðèì ïðÿìóþ $l$
ñ óãëîâûì êîýôôèöèåíòîì $t$, ïðîõîäÿùóþ ÷åðåç òî÷êè (0,-1)
è $(x,y):\quad y+1=tx.$

\begin{figure}
\begin{center}
\includegraphics[width=5cm,height=5cm]{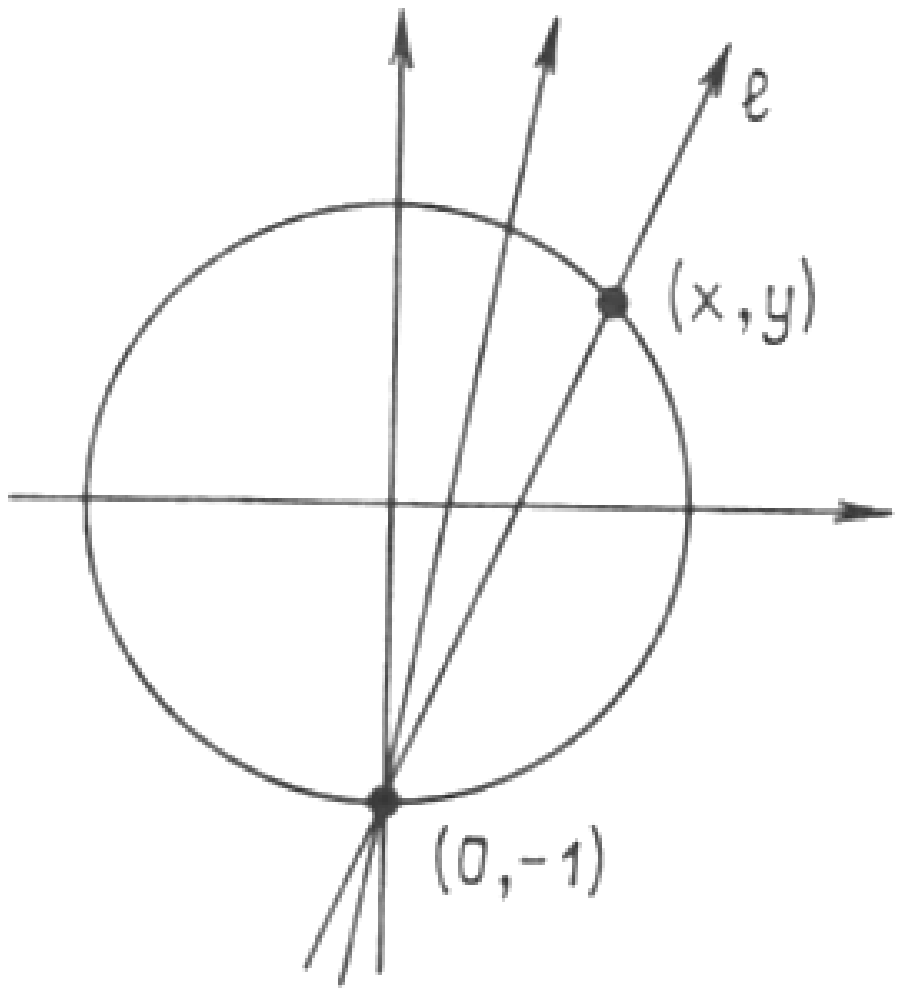}
\end{center}
\caption{\label{Img3}}
\end{figure}

Ïðè íàõîæäåíèè ðàöèîíàëüíûõ ðåøåíèé óðàâíåíèÿ
\begin{align}
 F(X_0,X_1,\dots ,X_n)=0 \label{i12.10}
\end{align}
 (ñ êâàäðàòè÷íîé ôîðìîé $F$ â (\ref{i12.8})) íàä $\Q$ 
 ìîæíî ñ÷èòàòü, ÷òî ôîðìà $F$ -- äèàãîíàëüíà:
ìåòîä Ëàãðàíæà âûäåëåíèÿ ïîëíûõ êâàäðàòîâ äàåò çàìåíó ïåðåìåííûõ $X=CY$
ñ ðàöèîíàëüíîé íåâûðîæäåííîé ìàòðèöåé $C\in \M_{n+1}(\Q)$.

Äëÿ îäíîðîäíûõ óðàâíåíèé òèïà (\ref{i12.10}) íåò ñóùåñòâåííîé ðàçíèöû ìåæäó
èõ öåëî÷èñëåííûìè è ðàöèîíàëüíûìè ðåøåíèÿìè: ïîñëå óìíîæåíèÿ íà
ïîäõîäÿùåå öåëîå ÷èñëî ëþáîå ðàöèîíàëüíîå ðåøåíèå ñòàíîâèòñÿ öåëî÷èñëåííûì, è åãî ìîæíî ñ÷èòàòü ïðèìèòèâíûì, ò. å. èìåþùèì âçàèìíî
ïðîñòûå â ñîâîêóïíîñòè êîîðäèíàòû. Íàèáîëåå ôóíäàìåíòàëüíûì ôàêòîì
òåîðèè êâàäðàòè÷íûõ ôîðì íàä ïîëåì ðàöèîíàëüíûõ ÷èñåë ÿâëÿåòñÿ ñëåäóþùèé ðåçóëüòàò.

\subsection{Ïðèíöèï ÌèíêîâñêîãîÕàññå äëÿ êâàäðàòè÷íûõ ôîðì.}
\label{i12.4.}
\index{Minkowski-Hasse principle!-- for a quadratic form}
\begin{theorem}
Öåëî÷èñëåííàÿ êâàäðàòè÷íàÿ ôîðìà $F(x_1,x_2,\dots ,x_n)$ ðàíãà $n$
ïðåäñòàâëÿåò íóëü íàä ïîëåì ðàöèîíàëüíûõ ÷èñåë òîãäà è òîëüêî
òîãäà, êîãäà äëÿ âñåõ íàòóðàëüíûõ ÷èñåë $N$, ñðàâíåíèå $F(x_1, \dots ,x_n)\equiv 0 \ (\bmod \  
N)$
èìååò ïðèìèòèâíîå
ðåøåíèå è ôîðìà
 $F$ ïðåäñòàâëÿåò
íóëü íàä ïîëåì âåùåñòâåííûõ ÷èñåë (ò. å. îíà íåîïðåäåëåííàÿ).
\end{theorem}
Ñì.
\cite{Borevich Z.I. Shafarevich I.R. (1985)}, ãëàâà 1. 
Êîíå÷íî, óòâåðæäåíèå «òîëüêî òîãäà» òðèâèàëüíî.

Ïðèâåäåì êðàñèâîå äîêàçàòåëüñòâî ýòîé òåîðåìû äëÿ ñëó÷àÿ, ðàññìîòðåííîãî Ëåæàíäðîì 
(\cite{Borevich Z.I. Shafarevich I.R. (1985)}):
Ïóñòü
$$
F=a_1x_1^2+a_2x_2^2+a_3x_3^2\qquad (a_1a_2a_3\neq 0).
$$
Íåîïðåäåëåííîñòü ôîðìû $F$ îçíà÷àåò, ÷òî íå âñå êîýôôèöèåíòû $F$ îäíîãî çíàêà. 
Óìíîæèâ ôîðìó ïðè íåîáõîäèìîñòè íà $-1$  ìû ïðèäåì ê ñëó÷àþ, êîãäà äâà êîýôôèöèåíòà ïîëîæèòåëüíû, à îäèí îòðèöàòåëåí.
Êðîìå
òîãî, ìû ìîæåì ñ÷èòàòü ýòè  ÷èñëà
 öåëûìè, ñâîáîäíûìè îò êâàäðàòîâ
è âçàèìíî ïðîñòûìè â ñîâîêóïíîñòè, òàê êàê èõ ìîæíî ñîêðàòèòü íà íàèáîëüøèé îáùèé äåëèòåëü. 
Äàëåå, åñëè, íàïðèìåð,
$a_1$ è $a_2$ èìåþò îáùèé
ïðîñòîé äåëèòåëü $p$, òî, óìíîæèâ ôîðìó íà $p$ è âçÿâ $px$ è $py$ çà íîâûå
ïåðåìåííûå, ìû ïîëó÷èì ôîðìó ñ êîýôôèöèåíòàìè $a_1/p$, $a_2/p$ è $pa_3$. 
Ïîâòîðÿÿ ýòîò ïðîöåññ íåñêîëüêî ðàç, ìû çàìåíèì íàøó ôîðìó ôîðìîé âèäà 
\begin{align}
 ax^2+by^2-cz^2. \label{i12.11}
\end{align}
â êîòîðîé öåëûå ïîëîæèòåëüíûå ÷èñëà $a,b, c$ ïîïàðíî âçàèìíî ïðîñòû
(è ñâîáîäíû îò êâàäðàòîâ). Ïóñòü òåïåðü
 $p$ êàêîé-íèáóäü ïðîñòîé äåëèòåëü ÷èñëà $c$, îòëè÷íûé îò $2$. 

Ìîæíî ïîêàçàòü, ÷òî ïîñêîëüêó äëÿ èñõîäíîé
ôîðìû ñóùåñòâóåò ïðèìèòèâíîå ðåøåíèå ñðàâíåíèÿ $F\equiv 0 (\bmod \  p^l)$ äëÿ ëþáîãî $l\ge 1$, 
 òî ñðàâíåíèå $ax^2+by^2\equiv 0 (\bmod \  p)$ 
èìååò íåòðèâèàëüíîå ðåøåíèå
$(x_0, y_0)$. 
Ñëåäîâàòåëüíî, ìîæíî ïðåäïîëàãàòü, ÷òî $y_0\ne 0$,
è âûïîëíÿåòñÿ ðàçëîæåíèå íà ìíîæèòåëè
$$
 ax^2+by^2\equiv ay_0^{-2}(xy_0+yx_0)(xy_0-yx_0)\ (\bmod \  p).
$$
Àíàëîãè÷íûå ðàçëîæåíèÿ èìåþò ìåñòî ïî ìîäóëþ íå÷åòíûõ $p$, äåëÿùèõ
êîýôôèöèåíòû $a$ è $b$, à ïðè $p=2$ âûïîëíÿåòñÿ ñðàâíåíèå
$$
ax^2+by^2-cz^2\equiv (ax+by-cz)^2\ (\bmod \  2).
$$
Òàêèì îáðàçîì, äëÿ ëþáîãî ïðîñòîãî ÷èñëà $p|2abc$ ñóùåñòâóþò ëèíåéíûå ôîðìû 
$L^{(p)},\ M^{(p)}$
îò  $x,y,z$ ñ öåëûìè êîýôôèöèåíòàìè, òàêèå, ÷òî $F\equiv 
L^{(p)}M^{(p)}
(\bmod \  p)$. Òåïåðü ñ ïîìîùüþ êèòàéñêîé òåîðåìû îá îñòàòêàõ
íàéäåì òàêèå ëèíåéíûå ôîðìû $L$ 
(ñîîòâ. $M$)
ñ öåëûìè êîýôôè\-öèåí\-òà\-ìè, ñðàâíèìûìè ñ  $L^{(p)}$ (ñîîòâ. 
$M^{(p)}$)
$(\bmod \  p)$ äëÿ âñåõ $p|abc$,
è ìû ïîëó÷èì
\begin{align}
ax^2+by^2+cz^2\equiv L(x,y,z)M(x,y,z) \ (\bmod \ abc). \label{i12.12}
\end{align}
Áóäåì ïðèäàâàòü ïåðåìåííûì $x$, $y$, $z$
 öåëûå çíà÷åíèÿ, óäîâëåòâîðÿþùèå
 óñëîâèÿì
\begin{align}
0\le x<\sqrt{bc},\quad 0\le y<\sqrt{ac},\quad 0\le z<\sqrt{ab}.\label{i12.13}
\end{align}
Åñëè èñêëþ÷èòü èç ðàññìîòðåíèÿ òðèâèàëüíûé ñëó÷àé $a=b=c=1$, 
òî íå âñå ÷èñëà $\sqrt{bc}, \sqrt{ac}, \sqrt{ab}$ öåëûå è ÷èñëî òðîåê $(x, y, z)$,
óäîâëåòâîðÿþùèõ óñëîâèÿì (\ref{i12.13}), ñòðîãî áîëüøå ÷åì îáú¸ì $\sqrt{bc} \sqrt{ac} \sqrt{ab} = abc$. 
Ñëåäîâàòåëüíî, äëÿ íåêîòîðûõ äâóõ ðàçëè÷íûõ òðîåê 
 ëèíåéíàÿ ôîðìà $L$ ïðèíèìàåò îäíî è òî æå çíà÷åíèå  $\bmod \ abc$,
  îòêóäà â ñèëó ëèíåéíîñòè ôîðìû èìååì
\begin{align}
L(x_0,y_0,z_0)\equiv 0 \ (\bmod \ abc)  \label{i12.14}
\end{align}
äëÿ íåêîòîðîãî ðåøåíèÿ $|x_0|\le \sqrt{bc},\quad |y_0|\le \sqrt{ac},\quad |z_0|\le
\sqrt{ab}$. 
Ïîýòîìó
\begin{align}
ax_0^2+by_0^2-cz_0^2\equiv 0 \ (\bmod \ abc) \label{i12.15}
\end{align}
  è èìåþò ìåñòî íåðàâåíñòâà
$$-abc < ax_0^2+by_0^2-cz_0^2 < 2abc.$$
Òàêèì îáðàçîì,  ñïðàâåäëèâî îäíî èç
äâóõ ðàâåíñòâ
\begin{align}
ax_0^2+by_0^2-cz_0^2=0  \label{i12.16}
\end{align}
èëè æå
\begin{align}
ax_0^2+by_0^2-cz_0^2=abc. \label{i12.17}
\end{align}
Â ñëó÷àå (\ref{i12.16}) òåîðåìà äîêàçàíà; 
åñëè æå âûïîëíåíî ðàâåíñòâî (\ref{i12.17}),
òî äîêàçàòåëüñòâî ñëåäóåò èç òîæäåñòâåííîãî ïðåîáðàçîâàíèÿ
$$
a(x_0z_0+by_0)^2+b(y_0z_0-ax_0)^2-c(z_0^2+ab)^2=0.
$$
Â ôîðìóëèðîâêå Ëåæàíäðà äèîôàíòîâî óðàâíåíèå
 $ax^2+by^2-cz^2=0$ èìååò íåòðèâèàëüíîå öåëî÷èñëåííîå ðåøåíèå â òîì è òîëüêî â òîì ñëó÷àå,
êîãäà êëàññû âû÷åòîâ
$$
bc \ (\bmod\,   a),\quad ac \ 
(\bmod \,  b),
\quad -ab \ (\bmod  \, c)
$$
 ÿâëÿþòñÿ êâàäðàòàìè.

Ìîæíî äîêàçàòü, ÷òî ðàöèîíàëüíàÿ êâàäðàòè÷íàÿ ôîðìà ðàíãà $\ge 5$
âñåãäà ïðåäñòàâëÿåò íóëü íàä $\Q$. 

Â îáùåì ñëó÷àå ñóùåñòâóþò ýôôåêòèâíûå ìåòîäû
(îñíîâàííûå íà ïðèíöèïå Ìèíêîâñêîãî -- Õàññå äëÿ êâàäðàòè÷íûõ ôîðì), ÷òîáû óñòàíîâèòü, ïðåäñòàâëÿåò ëè êâàäðàòè÷íàÿ ôîðìà ðàöèîíàëüíûé íóëü. 
Ýòè ìåòîäû îñíîâàíû íà èíôîðìàöèè,
êîòîðóþ ìîæíî èçâëå÷ü èç âåùåñòâåííûõ è êîíãðóýíöèàëüíûõ ðàññìîòðåíèé,
è èñïîëüçóþò ñèìâîë Ãèëüáåðòà.

\subsection{Ñèìâîë Ãèëüáåðòà.}\label{ii23.3.}
Â ýòîì ïóíêòå ìû äîïóñêàåì çíà÷åíèå $p=\infty$,
 è ñ÷èòàåì òîãäà, ÷òî $\Q_\infty=\R$.
Ñèìâîë Ãèëüáåðòà (ñèìâîë íîðìåííîãî
âû÷åòà)
 $$
 (a,b)
 = {{a,b}\choose p}
 = \left ( \frac {a,b} p \right )
 = (a,b)_p
 $$
 äëÿ $a,b \in \Q_p^\times$ îïðåäåëÿåòñÿ ðàâåíñòâîì 
 $$
 (a,b)
 =
 \begin{cases}
 1 , & \mbox{åñëè ôîðìà }ax^2 + by^2 - z^2 \mbox{
       èìååò}
       \\ &
\mbox{íåòðèâèàëüíîå ðåøåíèå â} \  \Q_p; \cr
 -1, & \mbox{â ïðîòèâíîì ñëó÷àå.}
 \end{cases}
 $$
ßñíî, ÷òî $(a,b)$
 çàâèñèò òîëüêî îò $a$ è $b$ ïî ìîäóëþ êâàäðàòîâ â $\Q_p$. Ñóùåñòâóåò
íåñèììåòðè÷íàÿ ôîðìà ýòîãî îïðåäåëåíèÿ. Èìåííî, $(a,b)=1$ òîãäà è òîëüêî òîãäà, êîãäà
\begin{align}\label{ii23.7}
 a=z^2-by^2 \hbox{ äëÿ íåêîòîðûõ } y,\ z\in \Q_p.
\end{align}
Äåéñòâèòåëüíî, èç ñîîòíîøåíèÿ (\ref{ii23.7}) ñëåäóåò, ÷òî $(1,y,z)$ íåòðèâèàëüíûé íóëü
 êâàäðàòè÷íîé ôîðìû $ax^2 + by^2 - z^2$.
Íàîáîðîò, åñëè, $(x_0, y_0, z_0)$ íåêîòîðûé íåòðèâèàëüíûé íóëü,
 òî îñòàëüíûå íóëè ïîëó÷àþòñÿ ñ ïîìîùüþ
 ãåîìåòðè÷åñêîãî ïðèåìà ïðîâåäåíèÿ ñåêóùèõ ÷åðåç òî÷êó $(x_0, y_0, z_0)$ 
 ñ íàïðàâëÿþùèì âåêòîðîì, èìåþùèì êîîðäèíàòû èç $\Q_p$. 
Ïîýòîìó ìîæíî ñ÷èòàòü, ÷òî $x_0\not=0$.
Ïîýòîìó $(y_0/x_0,z_0/x_0)$ óäîâëåòâîðÿþò ñîîòíîøåíèþ (\ref{ii23.7}).
\medskip

{\it Ëîêàëüíûå ñâîéñòâà ñèìâîëà Ãèëüáåðòà:}
%
\begin{eqnarray} & & 
\mbox{à)} \hskip1cm
\label{ii23.8a}
 (a,b) = (b,a);
 \\ & & \mbox{á)} 
\label{ii23.8b}
\hskip1cm
( a_1a_2, b) = (a_1, b)(a_2, b),\ \ (a, b_1b_2) = (a,b_1)(a, b_2); 
 \\ & & \mbox{â)} 
\label{ii23.8c}
 \hskip1cm
 \hbox{åñëè } (a,b)=1 
\hbox{ äëÿ âñåõ } b, \hbox{ òî } a \in 
\Q_p^\times{}^2;
 \\ & & \mbox{ã)} 
\label{ii23.8d}
\hskip1cm
 (a, -a) = 1 \hbox{ äëÿ âñåõ}\ a ;
 \\ & & \mbox{ä)} 
\label{ii23.8e}
\hskip1cm 
\hbox{ åñëè}\ p\not =2,\infty  \hbox{ è } |a|_p = |b|_p = 1,  
\hbox{ òî } (a,b)=1. 
\end{eqnarray}
\medskip
Â ÷àñòíîñòè, ïðè ôèêñèðîâàííîì $b$,
 âñå $a$, äëÿ êîòîðûõ $(a,b)=1$ îáðàçóþò ãðóïïó ïî óìíîæåíèþ. 
 Óðàâíåíèå (\ref{ii23.7}) âûðàæàåò òîò ôàêò, ÷òî $a$
ÿâëÿåòñÿ íîðìîé èç êâàäðàòè÷íîãî ðàñøèðåíèÿ 
 $\Q_p(\sqrt b)/\Q_p$
 (cf. 
\cite{Borevich Z.I. Shafarevich I.R. (1985)}, 
\cite{Serre J.--P. (1970)}).

Âû÷èñëåíèå ñèìâîëà Ãèëüáåðòà ïîçâîëÿåò ïîëíîñòüþ ðåøèòü «ãëîáàëüíûé» âîïðîñ î ïðåä\-ñòàâ\-ëåíèè íóëÿ ðàöèîíàëüíûìè êâàäðàòè÷íûìè
ôîðìàìè (ñ ïîìîùüþ òåîðåìû Ìèíêîâñêîãî -- Õàññå). Åñëè, ñêàæåì,
\begin{align}\label{ii23.9}
Q(x,y,z) = ax^2 + by^2 + cz^2\ \ (a,b,c\in \Q,\ c\not=0), 
\end{align}
òî ôîðìà (\ref{ii23.9}) ïðåäñòàâëÿåò íóëü íàä ïîëåì $\Q$ òîãäà è òîëüêî òîãäà, êîãäà âûïîëíÿåòñÿ ðàâåíñòâî 
$(-a/c, -b/c)_p = 1$ 
äëÿ âñåõ $p$ (âêëþ÷àÿ $p=\infty$).
1. Ýòîò êðèòåðèé ÿâëÿåòñÿ âåñüìà ýôôåêòèâíûì, òàê êàê $|a|_p = |b|_p = 1$ äëÿ ïî÷òè âñåõ $p$ ïðè÷åì â ýòîì ñëó÷àå $(a,b)_p=1$ åñëè $p\not=2, \infty$ ñîãëàñíî ñâîéñòâó (\ref{ii23.8e}).
Âûïèøåì òåïåðü òàáëèöó äëÿ $(a,b)_p$: 
\begin{table}[h]
\begin{center}
\caption{ \label{Table10} Ñèìâîë Ãèëüáåðòà äëÿ $p>2$.
Çäåñü $v$ îáîçíà÷àåò òàêîå ÷èñëî $v \in \Z$,
 ÷òî   $\left ( \ds \frac v  p \right) = -1$,
  à $\varepsilon = 1$ åñëè $-1 \in \Q_p^\times{}^2$
 (ò. å. åñëè $p\equiv 1 \bmod\ 4$), à $\varepsilon = -1$ â ïðîòèâíîì ñëó÷àå \newline \ 
 }
\begin{tabular} {|l r|c|c|c|r|}
\hline
  & $a$ & \quad 1 \quad & \quad $v$ \quad & \quad $p$ \quad & \quad $pv$ \quad  \\
$b$ & & & & &\\ \hline
 1 & & $+1$ & $+1$ & $+1$ & $+1$  \\ \hline
 $v$ & & $+1$ & $+1$ & $-1$ & $-1$  \\ \hline
 $p$ & & $+1$ & $-1$ & $\varepsilon$ & $-\ep$ \\ \hline
 $pv$ & & $+1$ & $-1$ & $-\ep$ & $\ep$ \\
\hline
\end{tabular}
\end{center}
\end{table}
%

\medskip

{\it Ãëîáàëüíîå ñâîéñòâî ñèìâîëà Ãèëüáåðòà (ôîðìóëà ïðîèçâåäåíèÿ).}  
Ïóñòü $a, b \in \Q^\times$.
Òîãäà $(a,b)_p=1$ äëÿ ïî÷òè âñåõ $p$ è 
\begin{align}\label{ii23.10}
 \prod_{p \ {\rm including}\ \infty}(a,b)_p=1. 
\end{align}
Ôîðìóëà  (\ref{ii23.10}) ðàâíîñèëüíà êâàäðàòè÷íîìó çàêîíó âçàèìíîñòè. 
Äåéñòâèòåëüíî, ïî ñâîéñòâó (\ref{ii23.8e}) èìååì $|a|_p = |b|_p = 1$ äëÿ ïî÷òè âñåõ
 $p$, è â ýòîì ñëó÷àå $(a,b)_p=1$ for $p\not=2, \infty$ â ñèëó 
(\ref{ii23.8e}).
Îáîçíà÷èì
ëåâóþ ÷àñòü ðàâåíñòâà (\ref{ii23.10}) ÷åðåç $f(a,b)$.
Ïî ñâîéñòâàì(\ref{ii23.8b}) èìååì
\begin{align*}&
 f(a_1a_2, b) = f(a_1, b) f(a_2, b),\cr &
 f(a, b_1b_2) = f(a,b_1)  f(a,b_2),
\end{align*}
è ìîæíî ïðîâåðèòü, ÷òî $f(a,b) = 1$
êîãäà $a$ è $b$ ïðîáåãàþò ìíîæåñòâî îáðàçóþùèõ ãðóïïû $\Q^\times$:
 $
 -1,\ 2,\ -q \hbox{ íå÷åòíîå ïðîñòîå ÷èñëî}.
 $


\begin{table}
\caption{ \label{Table11} Ñèìâîë Ãèëüáåðòà â ñëó÷àå $p=2$.}
\begin{center}
\begin{tabular} {l r|c|c|c|c|c|c|c|r}\hline
  & $a$ & \quad 1 \quad & \quad $5$ \quad & \quad $-1$ \quad & \quad $-5$  \quad & \quad $2$ \quad & \quad $10$ \quad & \quad $-2$ \quad & \quad $-10$ \quad  \\
$b$ & & & & & & & & & \\ \hline
 1 & & $+1$ & $+1$ & $+1$ & $+1$  & $+1$ & $+1$ & $+1$ & $+1$  \\ \hline
 5 & & $+1$ & $+1$ & $+1$ & $+1$  & $-1$ & $-1$ & $-1$ & 
$-1$  \\ \hline
 $-1$ & & $+1$ & $+1$ & $-1$ & $-1$  & $+1$ & $+1$ & $-1$ & $-1$  \\ \hline
  $-5$ & & $+1$ & $+1$ & $-1$ & $-1$  & $-1$ & $-1$ & $+1$ & $+1$  \\ \hline
 $2$ & & $+1$ & $-1$ & $+1$ & $-1$  & $+1$ & $-1$ & $+1$ & 
$-1$  \\ \hline
 10 & & $+1$ & $-1$ & $+1$ & $-1$  & $-1$ & $+1$ & $-1$ & $+1$  \\ \hline
 $-2$  & & $+1$ & $-1$ & $-1$ & $+1$  & $+1$ & $-1$ & $-1$ & $+1$  \\ \hline
 $-10$ & & $+1$ & $-1$ & $-1$ & $+1$  & $-1$ & $+1$ & $+1$ & $-1$  \\ \hline
\end{tabular}

\end{center}
\end{table}

Îòìåòèì òàêæå ñëåäóþùåå ãëîáàëüíîå ñâîéñòâî
íîðìèðîâàíèé $|\cdot|_p$, àíàëîãè÷íîå ñâîéñòâó (\ref{ii23.10}):
\medskip

{\it Ôîðìóëà ïðîèçâåäåíèÿ äëÿ íîðìèðîâàíèé}.\index{Product formula 
for absolute values}
Ïóñòü $a\in \Q^\times$,
òîãäà $|a|_p=1$ äëÿ ïî÷òè âñåõ ïðîñòûõ ÷èñåë $p$,
 è
\begin{align}\label{ii23.11}
\prod_{p\ {\rm including}\ \infty}|a|_p=1. 
\end{align}
Äåéñòâèòåëüíî, åñëè $a\in \Q^\times$, òî 
$$
a = \pm \prod_{p\not=\infty}p^{v_p(a)},
$$
ãäå $v_p(a) \in \Z$ è $v_p(a)$ äëÿ ïî÷òè âñåõ $p$.
Òîãäà 
$$
|a|_p = p^{-v_p(a)} \hbox{ (äëÿ $p\not=\infty$),}
$$
$$
|a|_\infty = \prod_{p\not=\infty} p^{v_p(a)}.
$$

\section {Êóáè÷åñêèå óðàâíåíèÿ è ýëëèïòè÷åñêèå êðèâûå}\label{i13.}

\subsection{Ïðîáëåìà ñóùåñòâîâàíèÿ ðàöèîíàëüíîãî ðåøåíèÿ.}\label{i13.1.}
Äëÿ öåëî÷èñëåííûõ êóáè÷åñêèõ ôîðì
 $F(X,Y,Z)$ îò òðåõ ïåðåìåííûõ óæå íå
èçâåñòíî íèêàêîãî îáùåãî àëãîðèòìà, ïîçâîëÿþùåãî óñòàíîâèòü ñóùåñòâîâàíèå íåòðèâèàëüíîãî ðåøåíèÿ íàä $\Q$,
õîòÿ èçó÷åíî áîëüøîå ÷èñëî
êîíêðåòíûõ óðàâíåíèé, íàïðèìåð óðàâíåíèé âèäà
$$
 aX^3+bY^3+cZ^3=0.
$$
Îêàçûâàåòñÿ, äëÿ êóáè÷åñêèõ ôîðì ïåðåñòàåò,
âîîáùå ãîâîðÿ, âûïîëíÿòüñÿ ïðèíöèï Ìèí\-êîâñ\-êîãî -- Õàññå: óðàâíåíèå
 $3X^3+4Y^3+5Z^3=0$
 íå èìååò íåòðèâèàëüíûõ ðåøåíèé â öåëûõ ÷èñëàõ,
õîòÿ èìååò âåùåñòâåííûå ðåøåíèÿ, è äëÿ âñåõ íàòóðàëüíûõ $N>1$ ñðàâíåíèå
$ 3X^3+4Y^3+5Z^3=0\bmod N$ èìååò ïðèìèòèâíîå ðåøåíèå. 
Íàðóøåíèå ïðèíöèïà ÌèíêîâñêîãîÕàññå ìîæåò áûòü èçìåðåíî ÷èñëåííî ïðè
ïîìîùè ãðóïïû {\it Øàôàðåâè÷àÒýéòà}, ñì. ãëàâó 5 êíèãè \cite{Ma-Pa05}.

\subsection{Ñëîæåíèå òî÷åê íà êóáè÷åñêîé êðèâîé.}\label{i13.2.}
\index{Group!-- of rational points on a non-singular cubic curve}
\index{Curve!Cubic --}
Êóáè÷åñêàÿ ôîðìà
 $F(\! X,\! Y,\! Z)$ çàäàåò êðèâóþ ${\cal  C}$
 íà ïðîåêòèâíîé ïëîñêîñòè ${\Bbb P}^2$:
\begin{align}\label{i13.1}
{\cal  C} = \{ (X:Y:Z)\ |\ F(X,Y,Z)=0\} . 
\end{align}
ïðè÷åì ìû ñ÷èòàåì, ÷òî êîîðäèíàòû â ôîðìå (\ref{i13.1})
êîìïëåêñíûå ÷èñëà. Åñëè íà  ${\cal  C}$ ëåæèò õîòÿ áû îäíà ðàöèîíàëüíàÿ òî÷êà ${O}$, è  êðèâàÿ ${\cal  C}$  íåâûðîæäåíà, òî ìîæíî íàéòè òàêóþ îáðàòèìóþ çàìåíó êîîðäèíàò (íàä
ïîëåì $\Q$) ïîñëå êîòîðîé ôîðìà $F$ ïðèìåò âèä
\begin{align}\label{i13.2}
 Y^2Z-X^3-aXZ^2-bZ^3 \qquad (a,b \in {\Bbb Q}).
\end{align}
(âåéåðøòðàññîâà ôîðìà), ïðè÷åì òî÷êà ${O}$ ïåðåéäåò â ðåøåíèå $(0:1:0)$ 
äëÿ ôîðìû (\ref{i13.2}), à óñëîâèå íåâûðîæäåííîñòè äëÿ êðèâîé (\ref{i13.2})
ñòàíåò ýêâèâàëåíòíî òîìó, ÷òî $4a^3+27b^2\ne 0$. 

Íåâûðîæäåííàÿ êóáè÷åñêàÿ êðèâàÿ, èìåþùàÿ ðàöèîíàëüíóþ òî÷êó, íàçûâàåòñÿ ýëëèïòè÷åñêîé êðèâîé. 
Â íåîäíîðîäíûõ êîîðäèíàòàõ $x=X/Z, y=Y/Z$
óðàâíåíèå êðèâîé $F=0$ ïðèìåò âèä
\begin{align}\label{i13.3}
 y^2=x^3+ax+b,
\end{align}
ïðè÷åì êóáè÷åñêèé ìíîãî÷ëåí ñïðàâà íå èìååò êðàòíûõ êîìïëåêñíûõ êîð-
íåé (åãî äèñêðèìèíàíò îòëè÷åí îò íóëÿ), à òî÷êà $O=(0 : 1 : 0)$
â ýòîé çàïèñè ñòàíåò áåñêîíå÷íî óäàëåííîé òî÷êîé. 
Ñóùåñòâóåò êðàñèâûé
ãåîìåòðè÷åñêèé ñïîñîá ïðåâðàòèòü ìíîæåñòâî ðàöèîíàëüíûõ òî÷åê
 ${\cal  C}$ íà òàêîé
êðèâîé â àáåëåâó ãðóïïó ñ íåéòðàëüíûì ýëåìåíòîì  $O$(«ìåòîä
ñåêóùèõ è êàñàòåëüíûõ»), ñì.
\cite{Shafarevich I.R. (1988)}, 
\cite{Ma-Pa05}.
Åñëè $P,Q\in {\cal  C}({\Bbb Q})$,
 òî ïðîâîäèì ÷åðåç  $P,Q$  ïðîåêòèâíóþ ïðÿìóþ, ïåðåñåêàþùóþ ${\cal  C}$ 
 â îäíîçíà÷íî îïðåäåëåííîé òðåòüåé òî÷êå $P'\in {\cal  C}({\Bbb Q})$.
çàòåì ïðîâîäèì ïðÿìóþ ÷åðåç $P^\prime$ è $O$ ,
à òî÷êó åå ïåðåñå÷åíèÿ ñ ${\cal  C}$ íàçîâåì ñóììîé òî÷åê $P+Q$.
Àíàëîãè÷íî îïðåäåëÿåòñÿ òî÷êà $2P$ åñëè èñïîëüçîâàòü êàñàòåëüíóþ, ïðîõîäÿùóþ ÷åðåç òî÷êó 
 $P.$ 
\begin{figure}[h]
\hskip1cm
\includegraphics[width=5cm,height=5cm]{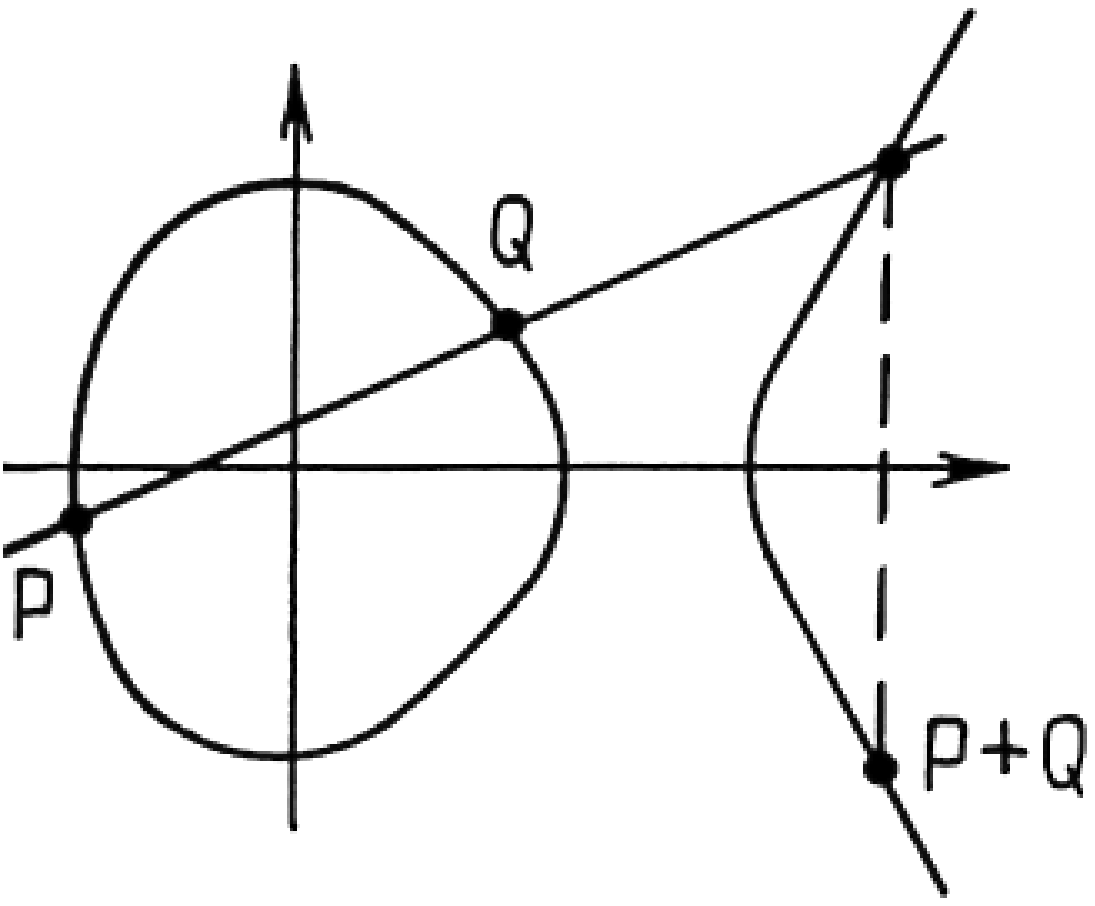}
\hskip2cm
\includegraphics[width=5cm,height=5cm]{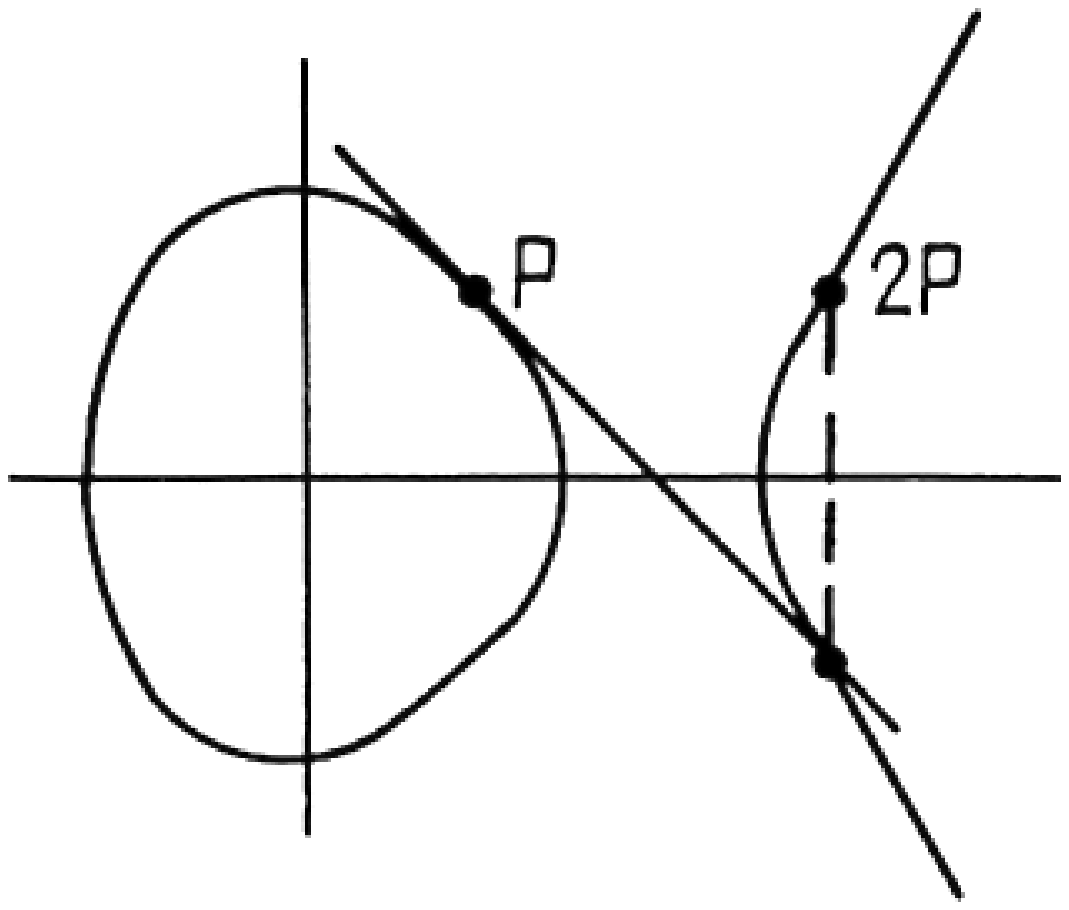}
\caption{\label{Img6}} 
\caption{\label{Img7}}
\end{figure}

Åñëè
 $P=(x_1,y_1),\ 
Q=(x_2,y_2)$ â íåîäíîðîäíûõ êîîðäèíàòàõ, ïðè÷åì $x_1\ne x_2$, 
òî $P+Q=(x_3,y_3)$, ãäå
\begin{align} \label{i13.4}
&x_3=-x_1-x_2+\left( \frac{y_1-y_2} {x_1-x_2} \right)^2, 
\cr
&y_3=\frac{y_1-y_2} {x_1-x_2} (x_1-x_3)-y_1. 
\end{align} 

Åñëè æå $P=Q$, òî
\begin{align}\label{i13.5}
x_3=-2x_1+\left( \frac{3x_1^2+a} {2y_1} \right)^2,
\ y_3=\frac{3x_1^2+a} {2y_1} (x_1-x_3)-y_1. 
\end{align}
Åñëè $x_1=x_2$, íî $y_1=-y_2$ òî òî÷êà $P+Q=O,$ áåñêîíå÷íî óäàëåííàÿ; îíà
âûáðàíà íåéòðàëüíûì ýëåìåíòîì ãðóïïîâîãî çàêîíà, ïîýòîìó â äàííîì
ñëó÷àå $P=-Q$.

Îïèñàííûé ìåòîä äàåò âîçìîæíîñòü ðàçìíîæàòü ðàöèîíàëüíûå òî÷êè, 
$mP,\ m\in {\Bbb Q}$
ðàññìàòðèâàÿ êðàòíûå  $mP,\ m\in {\Bbb Q}$, à òàêæå èõ ñóììû ñ äðóãèìè
òî÷êàìè $Q$ (åñëè òàêîâûå èìåþòñÿ).

Äëÿ âûðîæäåííûõ êóáè÷åñêèõ êðèâûõ îïèñàííûé ìåòîä íåïðèìåíèì. 
Ïóñòü,
ê ïðèìåðó,
\begin{align}\label{i13.6}
 {\cal  C} :\ y^2=x^2+x^3,
\end{align}

\begin{figure}
\begin{center}
\includegraphics[width=5cm,height=5cm]{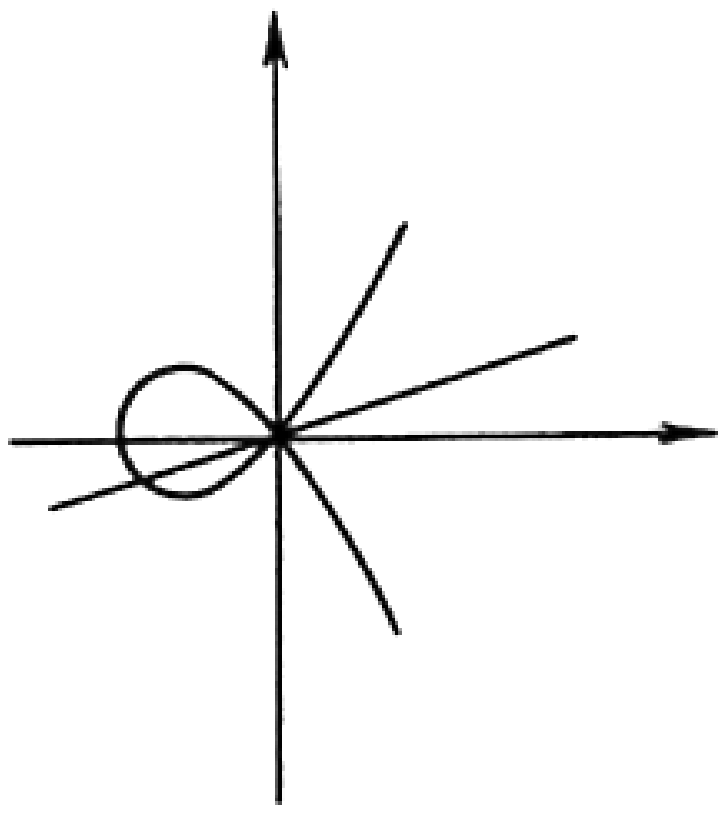}
\end{center}
\caption{\label{Img8}}

\end{figure}

\noindent
êðèâàÿ, èçîáðàæåííàÿ íà ðèñ. \ref{Img8}.
Òîãäà ëþáàÿ ïðÿìàÿ, ïðîõîäÿùàÿ ÷åðåç òî÷êó $(0,0)$ 
èìååò ëèøü îäíó òî÷êó ïåðåñå÷åíèÿ ñ êðèâîé
${\cal  C}$: åñëè óðàâíåíèå ïðÿìîé $y=tx$ òî ìû ïîëó÷àåì èç
óðàâíåíèÿ, ÷òî  $x^2(t^2-x-1)=0$.
Êîðåíü $t=0$ ñîîòâåòñòâóåò òî÷êå $(0, 0)$; $x=0,$ êðîìå òîãî, ìû èìååì åùå îäèí êîðåíü $x=t^2-1$.
 Èç óðàâíåíèÿ ïðÿìîé ìû ïîëó÷àåì, ÷òî 
$y=t(t^2-1)$
Ïîýòîìó, õîòÿ è íåëüçÿ îïðåäåëèòü ãðóïïîâîé çàêîí, êàê âûøå, ìû íàõîäèì âñå
ðàöèîíàëüíûå òî÷êè íà ${\cal  C}$ ñ ïîìîùüþ ðàöèîíàëüíîé ïàðàìåòðèçàöèè:$(x, y)=(t^2-1, t(t^2-1))$.

Âîîáùå, êðèâàÿ, äîïóñêàþùàÿ ïàðàìåòðèçàöèþ ñ ïîìîùüþ íåêîòîðûõ ðàöèîíàëüíûõ
ôóíêöèé ñ êîýôôèöèåíòàìè èç ïîëÿ $K$, íàçûâàåòñÿ ðàöèîíàëüíîé íàä $K$.\index{Curve!rational}

\subsection{Ñòðîåíèå ãðóïïû ðàöèîíàëüíûõ òî÷åê íà êóáè÷åñêîé êðèâîé}\label{i13.3.}
Íàèáîëåå âûäàþùàÿñÿ îñîáåííîñòü ìåòîäà ñåêóùèõ è êàñàòåëüíûõ  ýòî âîçìîæíîñòü ñâîäèòü íàõîæäåíèå âñåõ ðàöèîíàëüíûõ ðåøåíèé
êóáè÷åñêîãî óðàâíåíèÿ (\ref{i13.3}) ê íàõîæäåíèþ ëèøü êîíå÷íîãî èõ ÷èñëà.
Òî÷íåå, èìååò ìåñòî ñëåäóþùèé ðåçóëüòàò.

\index{Mordell's theorem}
\begin{theorem}[òåîðåìà Ìîðäåëëà]
Àáåëåâà ãðóïïà ${\cal  C}({\Bbb Q})$ êîíå÷íî
ïîðîæäåíà.
\end{theorem}
(ñì. \cite{Cassels J.W.S. (1966)}, 
è ïðèëîæåíèå
Þ. È. Ìàíèíà ê \cite{Mumford D. (1974)}).
Ñîãëàñíî òåîðåìå î ñòðîåíèè êîíå÷íî ïîðîæäåííûõ àáåëåâûõ ãðóïï
èìååòñÿ ðàçëîæåíèå
$$
 {\cal  C}({\Bbb Q}) \cong \Delta \times {\Bbb Z}^r
$$
ãäå $\Delta$  êîíå÷íàÿ ïîäãðóïïà âñåõ òî÷åê êðó÷åíèÿ,
è ${\Bbb Z}^r$ -- ïðÿìàÿ ñóììà áåñêîíå÷íûõ öèêëè÷åñêèõ ãðóïï; ÷èñëî
 $r$ íàçûâàþò ðàíãîì êðèâîé ${\cal  C}$ íàä $\Q$.

Î ãðóïïå êðó÷åíèÿ $\Delta$
óæå äàâíî áûëî êîå-÷òî èçâåñòíî. Òàê, Íàãåëëü è ïîçäíåå Ëóòö ïîëó÷èëè ñëåäóþùèé èíòåðåñíûé ðåçóëüòàò,
 äàþùèé îäíîâðåìåííî ìåòîä äëÿ ÿâíîãî îïðåäåëåíèÿ òî÷åê êðó÷åíèÿ
êîíêðåòíûõ êðèâûõ: åñëè $P = (x_P, y_P)$ðàöèîíàëüíàÿ òî÷êà êðó÷åíèÿ íà
êðèâîé, çàäàííîé óðàâíåíèåì $y^2 = x^3 + ax + b$, òî å¸ êîîðäèíàòû $x_P$ è $y_P$
ÿâëÿþòñÿ öåëûìè ÷èñëàìè, ïðè÷¸ì $y_P$ èëè ðàâíî  0, èëè $y_P^2$ ðàâåí êàêîìó - íèáóäü äåëèòåëþ äèñêðèìèíàíòà $D = -4aõ^3 - 27b^2$ äàííîé êðèâîé.

Á. Ìàçóð äîêàçàë â 1976 ã., ÷òî ïîäãðóïïà $\Delta$ êðó÷åíèÿ  íàä ${\Bbb Q}$
ìîæåò áûòü èçîìîðôíà ëèøü
îäíîé èç ïÿòíàäöàòè ãðóïï:
\begin{align}\label{i13.7}
{\Bbb Z}/m{\Bbb Z}\ (m\le 10, m=12),\ 
{\Bbb Z}/2{\Bbb Z}\times {\Bbb Z}/2n{\Bbb Z}\ (n\le 4), 
\end{align}
ïðè÷åì âñå âîçìîæíîñòè ðåàëèçóþòñÿ (ñì. \cite{stein}, ãëàâà 6).

\


Âû÷èñëåíèå ðàíãà $r$ îñòà¸òñÿ îòêðûòîé ïðîáëåìîé.

\begin{figure}
\begin{center}
\includegraphics[width=7cm,height=5cm]{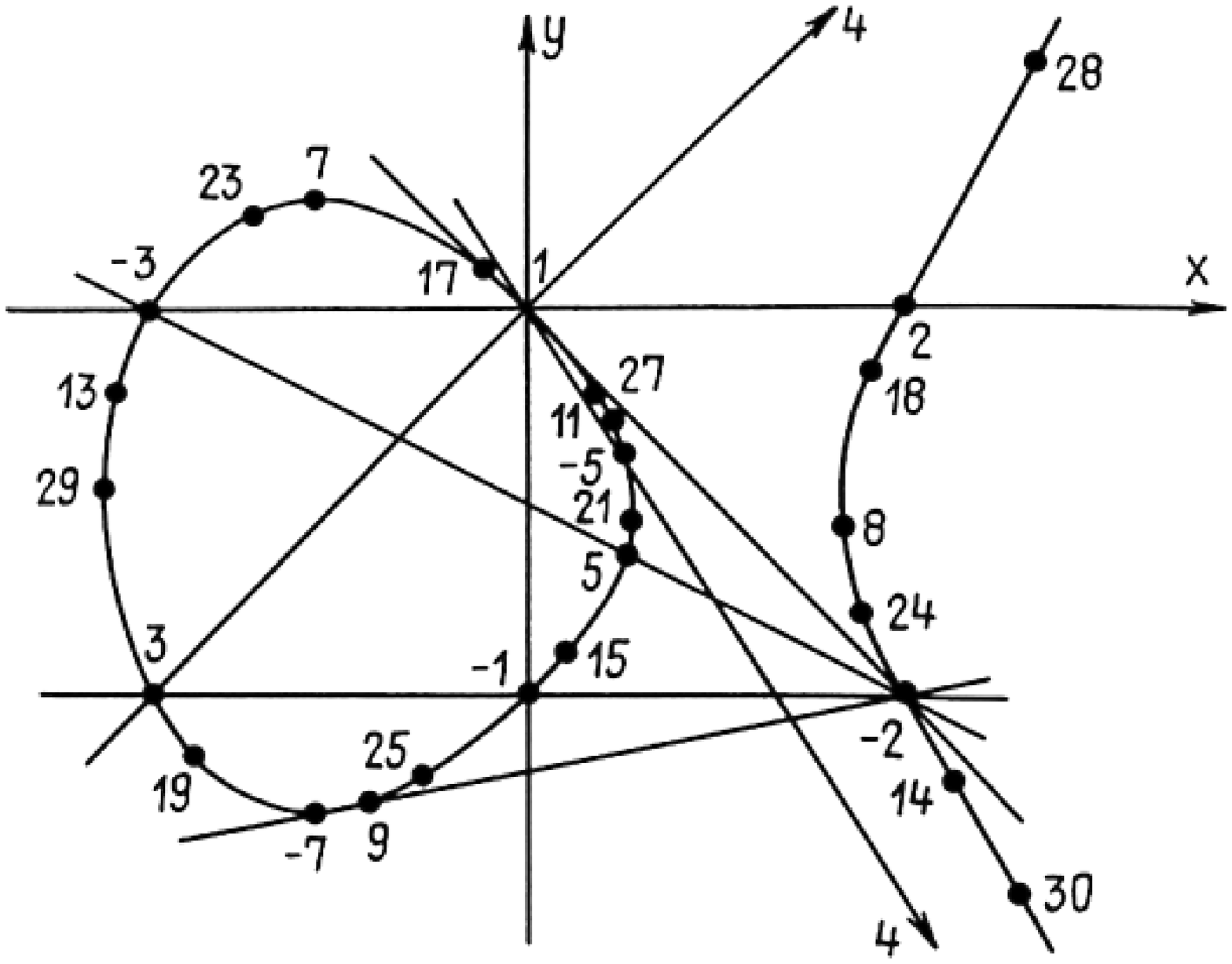}
\end{center}
\caption{\label{Img9}}
\end{figure}

\medskip
{\it Ïðèìåðû.}
1)
Êðèâàÿ ${\cal  C}$ çàäàåòñÿ óðàâíåíèåì 
$$
 y^2+y=x^3-x
$$
öåëî÷èñëåííîå ðåøåíèå êîòîðîãî äà¸ò ïðèìåð, êîãäà ïðîèçâåäåíèå äâóõ
ïîñëåäîâàòåëüíûõ öåëûõ ÷èñåë ðàâíî ïðîèçâåäåíèþ íåêîòîðûõ äðóãèõ òð¸õ
ïîñëåäîâàòåëüíûõ ÷èñåë. Òîãäà ãðóïïà $\Delta$ òðèâèàëüíà, à ãðóïïà òî÷åê ${\cal  C} ({\Bbb Q})$
(ñ áåñêîíå÷íî óäàëåííîé òî÷êîé â êà÷åñòâå íåéòðàëüíîãî ýëåìåíòà) ÿâëÿ-
åòñÿ áåñêîíå÷íîé öèêëè÷åñêîé ãðóïïîé, ïðè÷åì â êà÷åñòâå å¸ îáðàçóþùåé
ìîæíî âçÿòü òî÷êó $P=(0,0).$
Òî÷êè âèäà $mP$ óêàçàíû íà ðèñ.  \ref{Img9}.

2)
 Ïóñòü êðèâàÿ ${\cal  C}$ çàäàíà óðàâíåíèåì 
$$
y^2+y=x^3-7x+6.
$$
Òîãäà ${\cal  C}(\Q) \cong {\Bbb Z}^3$, à ñâîáîäíûå îáðàçóþùèå ýòîé ãðóïïû äàþòñÿ ðåøåíèÿìè (1,0), (2,0), 
(0, 2), ñì. \cite{BGZ85}.

3)
Êðèâàÿ ${\cal  C} : y(y+1)=x(x-1)(x+2)$ èìååò ðàíã, ðàâíûé äâóì, à êðèâàÿ $r=2;$ for ${\cal  C} : y(y+1)=x(x-1)(x+4),$ èìååò ðàíã $r=2$, ñð. ñ óðàâíåíèåì êðèâîé
ðàíãà 1 èç ïðèìåðà 1.

4) Ðàññìîòðèì êðèâóþ $y^2=x^3+px,\ p=877.$ Ìîæíî ïîêàçàòü (ñì. ññûëêè â), ÷òî  îáðàçóþùàÿ ïî ìîäóëþ êðó÷åíèÿ ãðóïïû ðàöèîíàëüíûõ òî÷åê íà ýòîé êðèâîé èìååò $x$-êîîðäèíàòó 
$$
x=\frac{375494528127162193105504069942092792346201} 
{6215987776871505425463220780697238044100} .
$$
Ýòîò ïðèìåð äàåò îïðåäåëåííîå ïðåäñòàâëåíèå î òðóäíîñòÿõ «íàèâíîãî» ýëåìåíòàðíîãî ïîäõîäà äëÿ íàõîæäåíèÿ òî÷åê áåñêîíå÷íîãî ïîðÿäêà íà ýëëèïòè÷åñêèõ êðèâûõ,ñì.
\cite{Coates J. (1984)}.

\subsection{Êóáè÷åñêèå ñðàâíåíèÿ ïî ïðîñòîìó ìîäóëþ.}\label{i13.4.}
\index{Curve!cubic --!reduction of -- modulo a prime}
Ïóñòü $p$ ïðîñòîå ÷èñëî, è 
 $F(X_0,X_1,$ $X_2)$ êóáè÷åñêàÿ ôîðìà, íåâûðîæäåííàÿ ïî ìîäóëþ $p$.
Ýòî çíà÷èò, ÷òî äëÿ ëþáîãî ïîëÿ $F\supset {\Bbb F}_p$
(ò. å. ïîëÿ õàðàêòåðèñòèêè $p$), ñëåäóþùèå ôîðìû ñòåïåíè 3 è 2:
$$
\overline F, {\partial \overline F\over \partial X_i} (i=0,1,2)
$$
íå èìåþò îáùèõ íåòðèâèàëüíûõ íóëåé íàä $K$, ãäå $\overline F$ 
îáîçíà÷àåò ôîðìó, ïîëó÷åííóþ èç $F$ ðàññìîòðåíèåì åå êîýôôèöèåíòîâ ïî ìîäóëþ $p$.

Êàê íàä ïîëåì ðàöèîíàëüíûõ ÷èñåë, ïðîñòûå àëãåáðî-ãåîìåòðè÷åñêèå
èäåè ìîæíî ïðèìåíèòü è ê ïîëþ $K$ ïîëîæèòåëüíîé õàðàêòåðèñòèêè. Â ýòîì
ñëó÷àå íîðìàëüíàÿ ôîðìà ñòàíîâèòñÿ íåñêîëüêî áîëåå ñëîæíîé. Ñäåëàâ
çàìåíó ïðîåêòèâíûõ êîîðäèíàò è ïåðåéäÿ ê íåîäíîðîäíîé ôîðìå çàïèñè,
ìû âñåãäà ìîæåì ïðèâåñòè óðàâíåíèå $F = 0$ ê âèäó
$$
  y^2 + a_1 xy + a_3 y = x^3 + a_2 x^2 + a_4 x + a_6,
$$
ãäå
$ a_1, a_2, a_3, a_4, a_6 \in K$
è
$$
  \Delta = -b_2^2 b_8 - 8b_4^3 - 27 b_6^2 + 9b_2 b_4 b_6 \neq 0,
$$
ãäå
$$
b_2 = a_1^2 + 4a_2, \quad b_4 = 2a_4 + a_1 a_3,\quad
  b_6 = a_3^2 + 4a_6.
$$
(Îáîçíà÷åíèÿ Òýéòà).
Òàêæå èñïîëüçóåòñÿ îáîçíà÷åíèå $\displaystyle j=\frac  {c_4^3}\Delta$, ãäå
$$
c_4=b_2^2-24b_4, c_6=-b_2^3+36b_2b_4-216b_6.
$$
Çàòåì ýòî óðàâíåíèå ìîæåò áûòü åùå óïðîùåíî ïðè ïîìîùè ïðåîáðàçîâàíèé âèäà
$x\mapsto u^2x'+r$,  $y\mapsto u^3y'+su^2x'r+t$ 
è ïîëó÷àåòñÿ ñëåäóþùåå (ñì.. 
\cite{Koblitz N. (1987)}: 
\begin{description}
\item{1)} Åñëè $p\ne 2,3$, òî
\begin{align}\label{i13.8}
y^2=x^3+a_4x+a_6\mbox{ ñ } \Delta=-16(4a_4^3+27a_6^2)\ne 0.
\end{align}
\item{2)} Åñëè $p=2$  óñëîâèå $j=0$ ðàâíîñèëüíî òîìó, ÷òî $a_1=0$,
è óðàâíåíèå
ïðåîáðàçóåòñÿ ê ñëåäóþùåìó âèäó: åñëè
 if $a_1\ne 0\  (i.e. j\ne 0)$, òî, âûáèðàÿ
ïîäõîäÿùèå $r,s,t$ ìû ìîæåì ïîëó÷èòü $a_1=1$, $a_3=0$, $a_4=0$, 
è óðàâíåíèå ïðèíèìàåò âèä
\begin{align}\label{i13.9'}
y^2+xy=x^3+a_2x^2+a_6,
\end{align} 
ãäå óñëîâèå ãëàäêîñòè çàäàåòñÿ ïðîñòî íåðàâåíñòâîì  $\Delta\ne 0$.
Ïðåäïîëîæèì òåïåðü, ÷òî $a_1= 0$  (ò.å. $j= 0$). Â ýòîì ñëó÷àå óðàâíåíèå ïðåîáðàçóåòñÿ â
\begin{align}\label{i13.9}
y^2+a_3y=x^3+a_4x+a_6,
\end{align}
è óñëîâèå ãëàäêîñòè â ýòîì ñëó÷àå çàäàåòñÿ íåðàâåíñòâîì $a_3\ne 0$.

\item{3)} Åñëè $p=3$, òî
\begin{align}\label{i13.10}
y^2=x^3+a_2x^2+a_4x+a_6,
\end{align}
\end{description}
(â ýòîì ñëó÷àå êðàòíûå êîðíè òàêæå íåäîïóñòèìû).
Â îäíîðîäíûõ êîîðäèíàòàõ âî âñåõ ñëó÷àÿõ äîáàâëÿåòñÿ «áåñêîíå÷íî óäàëåííîå ðåøåíèå» $O=(0:1:0).$

Êàê ïîñ÷èòàòü ÷èñëî ðåøåíèé òàêèõ 
 êóáè÷åñêèõ ñðàâíåíèé $F\equiv 0 \bmod \  p$? ßñíî, âî-ïåðâûõ, ÷òî ÷èñëî ýòèõ ðåøåíèé, îáðàçóþùèõ âìåñòå ñ òî÷êîé $O$ àáåëåâó ãðóïïó ïîðÿäîê êîòîðîé
íå ïðåâîñõîäèò $2p+1$,
òàê êàê äëÿ êàæäîãî òàêîãî $x$
íàéäóòñÿ íå áîëüøå äâóõ çíà÷åíèé $y$.
Îäíàêî ëèøü ïîëîâèíà ýëåìåíòîâ èç $(\F_p)^\times$ ÿâëÿþòñÿ êâàäðàòàìè, ïîýòîìó ìîæíî îæèäàòü, ÷òî ëèøü
â ïîëîâèíå ñëó÷àåâ èç ýëåìåíòà 
ìîæíî èçâëå÷ü êâàäðàòíûé
êîðåíü $y$ (ïðåäïîëîæèâ, ÷òî ýëåìåíòû
 $x^3+ax+b$ ðàçáðîñàíû ñëó÷àéíî â ïîëå $\F_p$. 

Áîëåå òî÷íî, ïóñòü
$\chi (x) = \left( \ds \frac {x} {p} \right)$
  ñèìâîë Ëåæàíäðà, îïðåäåëåíèå êîòîðîãî îçíà÷àåò, ÷òî ÷èñëî ðåøåíèé
óðàâíåíèÿ $y^2=u$ â ${\Bbb F}_p$
ðàâíî $1+\chi (u)$.
Òîãäà ìû ïîëó÷àåì ñëåäóþùóþ
ôîðìóëó äëÿ ÷èñëà ðåøåíèé êóáè÷åñêîãî ñðàâíåíèÿ:
\begin{align*}
 \Card \ {\cal  C} ({\Bbb F}_p)
 &=
 1+\sum_{x\in {\Bbb F}_p} (1+\chi (x^3+ax+b)) \cr
 &=
 p+1+\sum_{x\in {\Bbb F}_p} \chi (x^3+ax+b).
\end{align*}
Êîáëèö â \cite{Koblitz N. (1987)} ñðàâíèâàåò âçÿòèå ñóììû ñî «ñëó÷àéíûì áëóæäàíèåì», ïðè êîòîðîì äåëàåòñÿ øàã âïåðåä, åñëè $\chi (x^3+ax+b)=1$, è øàã
íàçàä, åñëè $\chi (x^3+ax+b)=-1$.
Èç òåîðèè âåðîÿòíîñòåé èçâåñòíî, ÷òî
ðàññòîÿíèå îò èñõîäíîé òî÷êè ïîñëå $p$  øàãîâ ïðè ñëó÷àéíîì áëóæäàíèè
áóäåò èìåòü ïîðÿäîê  $\sqrt{p}$.
È äåéñòâèòåëüíî, ýòî òàê: ñóììà âñåãäà îãðàíè÷åíà âåëè÷èíîé $2\sqrt{p}$.

\begin{theorem}[òåîðåìà Õàññå]\label{Hasse}
Ïóñòü $N_p=\Card\ {\cal  C }({\Bbb F}_p)$,
òîãäà
$$|N_p-(p+1)|\le 2\sqrt{p}.$$
\end{theorem}

Ýëåìåíòàðíîå äîêàçàòåëüñòâî ýòîãî ôàêòà áûëî äàíî 
Þ.È.Ìàíèíûì â  1956.

\subsection{Îò ñðàâíåíèé ê ðàöèîíàëüíûì òî÷êàì: ãèïîòåçà Á¸ð÷à è Ñóèííåð\-òîíà--Äàéåðà}
Çíàìåíèòûé ïðèìåð ñâÿçûâàþùèé ëîêàëüíóþ è ãëîáàëüíóþ èíôîðìàöèþ, äà¸òñÿ 
ãèïîòåçîé Á¸ð÷à è Ñóèííåðòîíà--Äàéåðà äëÿ ýëëèïòè÷åñêèõ êðèâûõ. 
Ýòà ãèïîòåçà ïðèíàäëåæèò ê ÷èñëó Ñåìè Ïðîáëåì Òûñÿ÷åëåòèÿ èíñòèòóòà CLAY,
à çà å¸ ðåøåíèå ïðåäëîæåí ïðèç â ìèëëèîí äîëëàðîâ!

Ýòà îòêðûòàÿ ïðîáëåìà îáñóæäàåòñÿ òàêæå â ñòàòüå Óàéëñà \cite{Wiles}.
\maketitle

\subsubsection{Ãèïîòåçà ÁÑÄ } 
(ñì. èçëîæåíèå â \cite{Stein}, ãëàâû 8 è 9)
Ïóñòü $E$ ýëëèïòè÷åñêàÿ êðèâàÿ íàä ~$\Q$ çàäàííàÿ óðàâíåíèåì
$$
  y^2 = x^3 + ax + b
$$
ñ $a,b\in\Z$.  Äëÿ $p\nmid \Delta = -16(4a^3 + 27b^2)$,
ïîëîæèì 
$a_p = p+1 - \# E(\Z/p\Z)$.
Ïóñòü 
\begin{align}\label{elll}
  L(E,s) = \prod_{p\nmid\Delta} \frac{1}{1-a_p p^{-s} + p^{1-2s}}
\end{align}
ðÿä Äèðèõëå, ñõîäÿùèéñÿ àáñîëþòíî ïðè $\Re (s) > \frac 3 2$ â ñèëó òåîðåìû Õàññå (òåîðåìà \ref{Hasse}).
\begin{theorem}[Áð¸é, Êîíðàä, Äàéàìîíä, Òýéëîð, Óàéëñ]\mbox{}\\
 Ôóíêöèÿ $L(E,s)$ ïðîëæàåòñÿ äî àíàëèòè÷åñêîé ôóíêöèè 
íà âñåé êîìïëåêñíîé ïðîñêîñòè  ${\CC}$.
\end{theorem}
(ñì. Breuil, Conrad, Diamond, Taylor, Wiles, à òàêæå ññûëêè â \cite{Ma-Pa05}).
\begin{conjecture}[Á¸ð÷à è Ñóèííåðòîíà--Äàéåðà]
Ðàçëîæåíèå Òýéëîðà ôóíêöèè  $L(E,s)$ â $s=1$ èìååò âèä
$$
  L(E,s) = c(s-1)^r + \text{\mbox{ ÷ëåíû âûñøåé ñòåïåíè}}
$$
ñ $c\neq 0$ è $E(\Q)\cong \Z^r \cross E(\Q)_{{\rm tors}}$.
\end{conjecture}

Ñïåöèàëüíûé ñëó÷àé ãèïîòåçû ÁÑÄ  óòâåðæäàåò, ÷òî $L(E,1)=0$ òîãäà è òîëüêî òîãäà, êîãäà $E(\Q)$ áåñêîíå÷íà,  
â ÷àñòíîñòè óòâåðæäåíèå ``$L(E,1)=0$ âëå÷¸ò, ÷òî ãðóïïà  $E(\Q)$ áåñêîíå÷íà''.

\subsubsection{×òî èçâåñòíî î ãèïîòåçå ÁÑÄ}
Â ñòàòüå \cite{Wiles} îáñóæäàåòñÿ èñòîðèÿ ñëåäóþùåãî ðåçóëüòàòà:
\begin{theorem}[Ãðîññ, Êîëûâàãèí, Çàãèð è äð.]
Ïðåäïîëîæèì, ÷òî 
$$
  L(E,s) = c(s-1)^r + \text{\mbox{ ÷ëåíû âûñøåé ñòåïåíè}}
$$
ñ  $r\leq 1$.  
Òîãäà ãèïîòåçà Á¸ð÷à è Ñóèííåðòîíà--Äàéåðà ñïðàâåäëèâà äëÿ~$E$, òî åñòü
$E(\Q)\cong \Z^r\oplus E(\Q)_{{\rm tors}}.$
\end{theorem}

\subsubsection{Âû÷èñëåíèÿ ñ ýëëèïòè÷åñêèìè êðèâûìè}
Îïèøåì, êàê  èñïîëüçîâàòü êîìïüþòåð  äëÿ ïðèáëèæ¸ííîãî âû÷èñëåíèÿ ðàíãà $r$ êðèâîé.

Ïóñòü~$E$ -- ýëëèïòè÷åñêàÿ êðèâàÿ íàä ïîëåì~$\Q$, îïðåä¸ëåííàÿ îáîáù¸ííûì óðàâíåíèåì Âåéåðøòðàññà
$$
  y^2 + a_1 xy + a_3 y = x^3 + a_2 x^2 + a_4 x + a_6.
$$
Íàïå÷àòàâ {\tt e = ellinit([$a_1$,$a_2$,$a_3$,$a_4$,$a_6$])}, ìû çàäàäèì ýòó êðèâóþ íà êîìïüþòåðíîé ñèñòåìå PARI
(ñì. \cite{BBBCO}).
Íàïðèìåð, íàïå÷àòàâ {\tt e = ellinit([0,0,1,-7,6])}, ìû ïîëó÷èì íà  PARI êðèâóþ
 $y^2 + y =x^3-7x+6$.
Ïðèâåä¸ì ïðèìåð âû÷èñëåíèÿ íà  PARI:

\subsubsection{Ñ ÷åãî íà÷àòü âû÷èñëåíèÿ íà PARI-GP}
\subsubsection*{Äîêóìåíòàöèÿ}
Äîêóìåíòàöèÿ äëÿ PARI äîñòóïíà ïî àäðåñó:
\begin{verbatim}
        http://pari.math.u-bordeaux.fr
\end{verbatim}
Âîò íåêîòîðàÿ äîêóìåíòàöèÿ äëÿ PARI:
\begin{enumerate}
\item {\bf Installation Guide:} Ïîìîùü ïî óñòàíîâêå  PARI íà êîìïüþòåðå.
\item {\bf Tutorial:} Ïðåâîñõîäíûé ââîäíûé òåêñò íà 42 ñòðàíèöû, êîòîðûé íà÷èíàåòñÿ ñ  {\tt 2 + 2}.
\item {\bf User's Guide:} Ïîäðîáíîå îïèñàíèå âñåõ ôóíêöèé íà 226 ñòðàíèöû.
\item {\bf Reference Card:} 4 ñòðàíèöû (î÷åíü ïîëåçíàÿ ñâîäêà êîìàíä è èõ ïðèìåíåíèÿ) 
\end{enumerate}
\begin{verbatim}
gp > factor(2^256+1)
%1=
[1238926361552897 1]

[93461639715357977769163558199606896584051237541638188580280321 1]

  ***   last result computed in 14,501 ms.
\end{verbatim}

\subsubsection{Êàê âû÷èñëèòü  $L(E,s)$ íà êîìïüþòåðå}
Ïóñòü~$E$ -- ýëëèïòè÷åñêàÿ êðèâàÿ íàä ïîëåì ~$\Q$, îïðåäåë¸ííàÿ îáîáù¸ííûì óðàâíåíèåì Âåéåðøòðàññà
$$
  y^2 + a_1 xy + a_3 y = x^3 + a_2 x^2 + a_4 x + a_6.
$$
Èìååòñÿ ìíîãî âîçìîæíîñòåé äëÿ âûáîðà óðàâíåíèÿ Âåéåðøòðàññà
îïðåäåëÿþøèõ ýëëèïòè÷åñêóþ êðèâóþ~$E$ ñ òî÷íîñòüþ äî èçîìîðôèçìà.
Ñðåäè ýòèõ óðàâíåíèé èìååòñÿ íàèëó÷øåå (ìèíèìàëüíîå), òî åñòü ñ íàèìåíüøèì äèñêðèìèíàíòîì.  
Ïðèìåð âû÷èñëåíèÿ íà  PARI:
\begin{verbatim}
? E = ellinit([0,0,0,-43,166]);
? E.disc
%61 = -6815744
? E = ellchangecurve(E,ellglobalred(E)[2])
%62 = [1, -1, 1, -3, 3, ...]
? E.disc
%63 = -1664
\end{verbatim}
Òàêèì îáðàçîì, óðàâíåíèå $y^2 + xy +y = x^3 -x^2 -3x +3$ ``ëó÷øå'', 
÷åì $y^2 = x^3 - 43x + 166$.  

{ÏÐÅÄÓÏÐÅÆÄÅÍÈÅ:} 
Íåêîòîðûå âàæíûå ôóíêöèè íà PARI äàþò ëèøü òîãäà âåðíûé ðåçóëüòàò, êîãäà óðàâíåíèå êðèâîé âûáðàíî
 íàèëó÷øèì (ìèíèìàëüíûì). Ýòî îòíîñèòñÿ ê òàêèì ôóíêöèÿì, êàê
{\tt elltors, ellap, ellak}, è {\tt elllseries}.

Íàïå÷àòàâ {\tt e = ellinit([0,0,1,-7,6])}, ìû ïîëó÷èì íà  PARI êðèâóþ
 $y^2 + y =x^3-7x+6$. Íàïå÷àòàâ  {\tt ellglobalred(e)} 
 ìû íàéä¸ì, ÷òî  óðàâíåíèå ìèíèìàëüíîå  è ÷òî êîíäóêòîð ðàâåí äèñêðèìèíàíòó 5077. 
\subsubsection{Ïðèáëèæ¸ííîå âû÷èñëåíèå ðàíãà}
Îïèøåì ìåòîä ïðèáëèæ¸ííîãî âû÷èñëåíèÿ ðàíãà êðèâîé íà  PARI.
Ìîæíî, íàïðèìåð, ïðèáëèæ¸ííî âû÷èñëÿòü 
çíà÷åíèÿ  $L^{(r)}(E,1)$ äëÿ $r=0,1,2,3,\dots$,
äî òåõ ïîð, ïîêà íå ïîëó÷èòñÿ íåíóëåâîå çíà÷åíèå. Ýòîò ìåòîä îïèñàí â êíèãå Äæ.Êðåìîíû  (ñì. \cite{Cremona}).

\begin{proposition}\label{propDLOG} 
Ïóñòü $L(E,s)=c(s-1)^r + \text{\mbox{ ÷ëåíû âûñøåé ñòåïåíè}}$.  
Òîãäà 
$$
\displaystyle{\lim_{s\rightarrow 1}}(s-1)\dfrac{L'(E,s)}{L(E,s)}=r\,.
$$
\end{proposition}

Òàêèì îáðàçîì, ðàíã $r$ âû÷èñëÿåòñÿ êàê ``òî÷íûé '' ïðåäåë ïðè  $s\rightarrow 1$ äëÿ ìåðîìîðôíûõ ôóíêöèé. 
Èçâåñòíî, ÷òî ýòîò ïðåäåë ÿâëÿåòñÿ öåëûì ÷èñëîì. Äëÿ êðèâîé $y^2 + y =x^3-7x+6$ 
ìîæíî ïîêàçàòü, ÷òî ýòîò ïðåäåë ðàâåí 3. Òåïåðü èñïîëüçóåì òàêîé ïðè¸ì: 
$$
(s-1)\frac{L'(s)}{L(s)}=\frac{s-1}{L(s)}\cdot\lim_{h\rightarrow 0}\frac{L(s+h)-L(s)}{h}\\
\approx \frac{s-1}{L(s)}\cdot\frac{L(s+(s-1)^2)-L(s)}{(s-1)^2}\\
=\frac{L(s^2-s+1)-L(s)}{(s-1)L(s)} .
$$

Ýòó ôîðìóëó âîçìîæíî èñïîëüçîâàòü íà  PARI äëÿ ïðèáëèæ¸ííîãî âû÷èñëåíèÿ ðàíãà $r$ êðèâîé.
{\small
\begin{verbatim}
gp > e=ellinit([0,0,1,-7,6]);
gp > r(E,s) = L1 = elllseries(E,s)
L2 = elllseries(E,s^2-s+1);
(L2-L1)/((s-1)*L1);
gp > r(e,1.00001)
%2 = 3.000011487248732705286325574
gp > ##
  ***   last result computed in 510 ms. 
  \end{verbatim}}
Íàïîìíèì, ÷òî ${\cal  C}(\Q) \cong {\Bbb Z}^3$, à ñâîáîäíûå îáðàçóþùèå ýòîé ãðóïïû äàþòñÿ ðåøåíèÿìè (1,0), (2,0), 
(0, 2), ñì. \cite{BGZ85}.

\

Äæîí Òýéò ñäåëàë äîêëàä î ãèïîòåçå ÁÑÄ  äëÿ èíñòèòóòà Clay. 
Ýòîò äîêëàä ìîæíî ïîñìîòðåòü îíëàéí ïî àäðåñó: 
\begin{center}
\hspace{-3.5em}
\begin{minipage}{1.1\textwidth}
\begin{verbatim}
http://www.msri.org/publications/ln/hosted/cmi/2000/cmiparis/index-tate.html
\end{verbatim}
\end{minipage}  
\end{center}

\subsection*{Ïðèçíàòåëüíîñòü àâòîðà}
Èñêðåííå áëàãîäàðþ Ýðíåñòà Áîðèñîâè÷à Âèíáåðãà  çà 
ïðèãëàøåíèå ïîäãîòîâèòü ñòàòüþ  äëÿ æóðíàëà ``Ìàòåìàòè÷åñêîå Ïðîñâåùåíèå'' 2008, 
ïîñâÿù¸ííîãî $p$-àäè÷åñêèì ÷èñëàì è èõ ïðè\-ëî\-æå\-íèÿì.


\bibliographystyle{plain}


\end{document}